\documentclass[11pt]{article}
\usepackage{amsmath,amsfonts,amssymb,bm}
\usepackage{tikz}

\usepackage{fancyhdr}
\usepackage{graphics}

\font\Ttl=cmb10 scaled \magstep2
\font\Srm=cmr9 
 
\font\Sbf=cmb9
\font\bfs=cmb10 scaled \magstep1

\def\RR{\mathbb{R}}

\def\QQ{\mathbf{Q}}

\font\Goth=eufm10
\font\goth = eufm6 
\def\SS{\hbox{{\Goth S}}}
\def\ss{\hbox{{\goth S}}}

\def\ppsi{\boldsymbol\psi}

\def\eps{\varepsilon}

\long\def\DEL#1{{}}

\begin{document}

\newcounter{SEC}
\newcounter{THM}
\def\ADS{\addtocounter{SEC}{1} \setcounter{THM}{0}}
\def\ADT{\addtocounter{THM}{1}}
\def\NUM{{\theSEC.\theTHM}}

%%%%%%%%%%%%%%%%%%%%%%%%%%%%%%%%%%%%u%%%%%%%%%%%%%%%%%%%%%%%%%%%%%%%%%%%%%%%%%%%%

\lhead{\textsc{L.L. Stach\'o}} 
\rhead{Locally generated $\mathcal{C}^1$-splines over triangular meshes}

\centerline{\Ttl Locally generated polynomial ${\cal C}^1$-splines} 

\centerline{\Ttl over triangular meshes}
\bigskip

\centerline{\textsc{L.L. STACH\'O}}
\bigskip
%\centerline{\it Communicated by $\phantom{
%\hbox{P. Hajnal}}$}

%REFERENCES

\def\Barn{1}
%R.E. Barnhill, Smooth interpolation over triangles (p.47) 45-70?, in:
%R.E. Barnhild and R.F. Riesenfeld (Ed.-s),
%Computer Aided Geometric Design, Academic Press, NY-San Francisco-London,1974

\def\Cox{2}
%H-S.M. Coxeter, Introduction to Geometry, New York: Wiley, 1969.

\def\Fulton{3}
%William Fulton (1974). Algebraic Curves. Mathematics Lecture Note Series. W.A. %Benjamin. p. 112. ISBN 0-8053-3081-4.

\def\Hamilton{4}
%H. Hamilton, Plane Algebraic Curves, Calderon Press, Oxford, 1920. 

\def\Lorenz{5} 
%Lorentz : Hermite interpolation by algebraic
%polynomials: A survey in Numerical Analysis 2000 (Elsevier)

\def\Meng{6}
%Meng wuMeng wuBernard MourrainBernard MourrainAndré GalligoB. NkongaB. %NkongaHermite Type Spline Spaces over Rectangular Meshes with Complex %Topological Structures
%March 2017Communications in Computational Physics 21(03):835-866
%DOI: 10.4208/cicp.OA-2016-0030

\def\Schumacher{7} 
%https://epubs.siam.org/doi/book/10.1137/1.9781611973907 
%L.L. Schumacher, Spline Functions, Computational Methods, SIAM 2015 

\def\Serg{8}
%I.V. Sergienko, O.M. Lytvyn, O.O. Lytvyn and O.I. Denisova,
%EXPLICIT FORMULAS FOR INTERPOLATING SPLINES OF DEGREE 5 ON THE TRIANGLE,
%Cybernetics and Systems Analysis, Vol. 50, No.5 (2014), 670-678.
%DOI 10.1007/s10559-014-9657-x

\def\Zlam{9} 
%M. Zl\'amal and A. \v Zeni\v sek, Mathematical aspects of the FEM, in
%Technical Physical and Mathematical Priciples of the FEM,
%Ed.-s V. Kola\v z et al., Akademia, Praha, 1971, pp. 15-39.

%%\def\St{10}
%%Home page http://www.math.u-szeged.hu/~stacho/

%%%%%%%%%%%%%%%%%%%%%%%%%%%%%%%%%%%%%%%%%%%%%%%%%%%%%%%%%%%%%%

\bigskip\bigskip
{\narrower \Srm \noindent
{\Sbf Abstract}. 
We classify all possible local linear procedures
over triangular meshes
resulting in polynomial $C^1$-spline functions 
with affinely uniform shape for the basic functions at the edges, 
and fitting the 9 value- and gradient data
at the vertices of the mesh members. 
There is a unique procedure among them with shape functions and
basic polynomials of degree 5 and all other admissible procedures are its pertubations
with higher degree.
\par}
\renewcommand{\thefootnote}{\fnsymbol{footnote}}
\footnote[0]{\hskip-6mm
Received March 29, 2019. \hfill\break
\hfill\break
}

\vskip10mm
\ADS
\noindent{\bfs \theSEC. Introduction}
\newcounter{Si} \setcounter{Si}{\value{SEC}}
\bigskip
\rm

By a {\it triangular mesh}
we mean a finite family
of closed non-degenerate triangles 
on the plane $\RR^2$ %(of real couples as usually) 
with pairwise non-intersecting interiors and
admitting only common vertices or edges. As usually,
we regard $\RR^2$ as the set of all real couples $[\xi,\eta]$
considered also as
$1\times 2$ (row) matrices. 
We shall use the standard notations
$x^{[1]}=x: [\xi,\eta] \mapsto \xi$, $x^{[2]}=y: [\xi,\eta] \mapsto \eta$
and $\langle {\bf u} \vert {\bf v} \rangle :=
\sum_{j=1}^2 x^{[j]}({\bf u}) x^{[j]}({\bf v})$ 
for the Cartesian
coordinates and scalar product, respectively. 
We write $\Vert {\bf u}\Vert =\langle {\bf u}\vert {\bf u}\rangle^{1/2}$
for the norm of ${\bf u}\in\RR^2$ and
${\rm Co}({\bf S})$ for the convex hull of 
${\bf S}\subset\RR^2$ resp. ${\rm det}({\bf u},{\bf v}) =
x^{[1]}({\bf u}) x^{[2]}({\bf v})-x^{[1]}({\bf v})x^{[2]}({\bf u})$  
for $2\!\times\! 2$-determinants.   
Given a triangular mesh
${\cal T} =\big\{ {\bf T}_1,\ldots,{\bf T}_M\big\}$,
in the sequel %we shall write
${\rm Vert}({\bf T}_k)$ and 
${\rm Edge}({\bf T}_k)$
will denote
%for 
the sets of vertices resp. closed edges
of the mesh members, 
furhermore \ 
${\rm Dom}({\cal T}) := \bigcup_{k=1}^M {\bf T}_k, \   
{\rm Edge}({\cal T}) := \bigcup_{k=1}^M \partial {\bf T}_k, \ 
% = \bigcup_{k=1}^M \Big[\bigcup {\bf E}_k \Big] , \  
{\rm Vert}({\cal T}) := \bigcup_{k=1}^M {\rm Vert}({\bf T}_k)$ 
will stand for the domain covered by 
${\cal T}$, 
the line figure covered by all edges and the collection of all vertices, respectively.
Recall that, given a {\it gradient-data}
%%\ADT
%%\newcounter{Ipi} \setcounter{Ipi}{\value{SEC}}
%%\newcounter{ipi} \setcounter{ipi}{\value{THM}}
$$F = \Big\{ \big( {\bf p}, f_{{\bf p}}, 
[ f^\prime_{x,{\bf p}},  f^\prime_{y,{\bf p}} ] \big) : {\bf p} 
\in {\rm Vert}{\cal T}) \Big\} \subset {\rm Vert}({\cal T})\times \RR\times\RR^2
$$   %%\eqno(\NUM)$$  %(1.1)$$
on the set of the vertices in ${\cal T}$,
a function $f: {\bf D} \to \RR$ is a
{\it ${\cal C}^1$-extension} of $F$ on ${\bf D}:={\rm Dom}({\cal T})$ if
$f$ has a continuous gradient 
${\bf p} \mapsto \nabla f({\bf p}) = 
\big[ \frac{\partial}{\partial x}f ({\bf p}), 
\frac{\partial}{\partial y}f  ({\bf p}) \big]$
on ${\rm Interior}({\bf D})$ which admits a
continuous extension to ${\bf D}$ as well 
(denoted also by $\nabla f$)
such that
\ADT
\newcounter{Ipii} \setcounter{Ipii}{\value{SEC}}
\newcounter{ipii} \setcounter{ipii}{\value{THM}} 
$$f({\bf p}) = f_{{\bf p}}, \quad 
\nabla f({\bf p}) = \big[ f^\prime_{x,{\bf p}} , f^\prime_{y,{\bf p}} \big]
\qquad 
\big({\bf p}\in{\rm Vert}({\cal T})\big) .
\eqno(\NUM)$$ %(1.2)$$
A ${\cal C}^1$-extension $f:{\bf D}\to\RR$ of $F$ 
is said to be a 
{\it $\mathcal{C}^1$-spline interpolation} 
of $F$ with respect to the mesh $\mathcal{T}$
if the restrictions $f\vert {\bf T}_k$ are polynomials 
of the coordinate functions $x,y$. 

There exists a large variety of ${\cal C}^1$-splines
for any admissible $\mathcal{T}$ and $F$ which can be obtained 
e.g. as global polynomial extensions with Hermite type interpolation [\Lorenz].
Obviously global polynomial fitting may primarily be interesting only from
a pure theretical view point 
due to too large polynomial degree and hence high numerical instability.
A better alternative could be an imitation of
tensor product splines 
(e.g. with Catmull-Rom type hermition curves on edges 
developed for rectangular meshes [\Schumacher,\Meng]). 
This consists the construction of 
$\mathcal{C}^1$-splines as linear combinations on the 
rectangular mesh members
from affine images of tensor products from only two special
polynomials $\Phi,\Psi:[0,1]\to[0,1]$ (actually
$\varphi(t)=t^2(3-2t)$, $\psi(t)=t^2(1-t)$).
Some main features of tensor product spline procedures 
which can naturally be generalized even to 
procedures
\ADT
\newcounter{Ipiii} \setcounter{Ipiii}{\value{SEC}}
\newcounter{ipiii} \setcounter{ipiii}{\value{THM}} 
$$\SS : ({\cal T},F) \!\mapsto\! f_{\mathcal{T},F} 
\quad \big(
\hbox{${\cal T}$ triang. mesh, $F$ 
grad. data on ${\rm Vert}({\cal T})$} \big)
\eqno(\NUM)$$ %(1.3)$$
furnishing ${\cal C}^1$-spline interpolation functions
from gradient data at the vertices 
over triangular meshes 
can be formulated in Postulates A,B below.

\medskip
\noindent
{\bf Postulate A.}
{\rm (Linearity and being locally generated)}. \ \it  
There are 
polynomial
functions
$$\varphi_{{\bf p},{\bf T}},\ 
\psi^{(1)}_{{\bf p},{\bf T}},\ \psi^{(2)}_{{\bf p},{\bf T}}: 
{\bf T}\!\to\!\RR
\quad \big(
{\bf T} \ \hbox{non-deg. triangle}\}, \  
{\bf p}\in {\rm Vert}({\bf T}) \big)$$
depending only on the couple of the 
triangle ${\bf T}$ with a distingvished vertex
such that 
the restriction of \ss\ \  to any mesh triangle ${\bf T}\in{\cal T}$ 
has the form
\ADT
\newcounter{Ipiv} \setcounter{Ipiv}{\value{SEC}}
\newcounter{ipiv} \setcounter{ipiv}{\value{THM}}
$$f_{{\cal T}\!,F} \vert {\bf T} = 
\sum_{{\bf p}\in {\rm Vert}({\bf T})} \Big[
f_{{\bf p}} \varphi_{{\bf p},{\bf T}} + 
f^\prime_{x,{\bf p}}  \psi^{(1)}_{{\bf p},{\bf T}} +
f^\prime_{y,{\bf p}}  \psi^{(2)}_{{\bf p},{\bf T}} \Big] .
\eqno(\NUM)$$ %(1.4)$$

%\noindent 
\rm
If Postulate A holds and
${\rm Vert}({\bf T}) =\{ {\bf a,b,p}\}$, 
in terms of the canonical frame vectors 
$${\bf e}^{[0]} :={\bf 0} =[0,0],\ {\bf e}^{[1]}:=[1,0],\ {\bf e}^{[2]}:=[0,1]$$
we necessarily have
\ADT
\newcounter{Ipv} \setcounter{Ipv}{\value{SEC}}
\newcounter{ipv} \setcounter{ipv}{\value{THM}}
$$\begin{aligned}
&\varphi_{{\bf p},{\bf T}} ({\bf p}) \!=\! 1, \ 
\nabla \varphi_{{\bf p},{\bf T}} ({\bf p}) \!=\! {\bf 0}, \qquad
\psi^{(j)}_{{\bf p},{\bf T}} ({\bf p}) \!=\! 0, \
\nabla \psi^{(j)}_{{\bf p},{\bf T}}({\bf p}) \!=\! {\bf e}^{[j]}; \\ 
&\varphi_{{\bf p},{\bf T}} ({\bf x}) \!=\!
\psi^{(j)}_{{\bf p},{\bf T}} ({\bf x}) \!=\! 0, \
\nabla\varphi_{{\bf p},{\bf T}} ({\bf x}) \!=\!
\nabla\psi^{(j)}_{{\bf p},{\bf T}} ({\bf x}) \!=\! {\bf 0}
\quad \big( {\bf x} \!\in\! {\rm Co}\{{\bf a, b}\} \big) .
\end{aligned}
\eqno(\NUM)$$ %(1.5)$$ 
The first statement in $(\theIpv.\theipv)$ %(1.5) 
is immediate from (\theIpiv.\theipiv), %(1.4), 
while
the second one is a consequence of the fact that given any point
$\widetilde{\bf p}$ forming an adjacent triangle 
$\widetilde{\bf T} := {\rm Co}\{ {\bf a,b},\widetilde{\bf p}\}$,
for the mesh ${\cal T} := \{ {\bf T},\widetilde{\bf T}\}$
with gradient data $F({\bf q}) = (0,{\bf 0})$ 
for ${\bf q}={\bf a,b},\widetilde{\bf p}$ we must have
$f_{{\cal T},F} \equiv 0$ on $\widetilde{\bf T}$ and hence
also $\nabla f_{{\cal T},F} \equiv 0$ on the common edge
${\rm Co}\{ {\bf a,b}\}$ of the triangles ${\bf T},\widetilde{\bf T}$.

Locally generated linear spline procedures have the computational advantage
that the resulting functions can be calculated on any mesh triangle 
regardless to what happens at vertices outside. 
A practical disadvantage is that in most cases only function values are
available (mostly from scanned data) and convenient gradient values 
must be guessed or found by optimizing procedres.

\medskip
{\bf Postulate B.}  
{\rm (Uniform shape on edges)}. \it 
$(\theIpiv.\theipiv)$ %$(1.4)$ 
holds and
there are polynomial functions \ 
$\Phi,\!\Psi:[0,1]\!\to\!\RR$ \ 
such that
\ADT
\newcounter{Ipvi} \setcounter{Ipvi}{\value{SEC}}
\newcounter{ipvi} \setcounter{ipvi}{\value{THM}} 
$$\Phi(0)\!=\!\Psi(0)\!=\!\Phi^\prime(0)\!=\!\Psi^\prime(0)\!=\!\Psi(1)=0, 
\quad 
\Phi(1)\!=\!\Psi^\prime(1)\!=\!1
\eqno(\theIpvi.\theipvi)$$ %(1.6)$$
and the graps of the basic functions 
$\varphi_{{\bf p},{\bf T}}$
on the edges of the triangle ${\bf T}$ 
are affine images of the graph of $\Phi$,
and those of $\psi^{(j)}_{{\bf p},{\bf T}}$ $(j=1,2)$
are affine images of the graph of $\Psi$.

\rm
That is, under Postulate B,
for the generic points ${\bf y}_t := t{\bf p}+(1-t){\bf a}$ 
on the edge
${\rm Co} \{ {\bf a},{\bf p}\}$,
resp. 
${\bf z}_t := t{\bf p}+(1-t){\bf b}$ on ${\rm Co} \{ {\bf b},{\bf p}\}$
we have 
\ADT
\newcounter{Ipvii} \setcounter{Ipvii}{\value{SEC}}
\newcounter{ipvii} \setcounter{ipvii}{\value{THM}}
$$\begin{aligned}
&\varphi_{{\bf p},{\rm Co}\{ {\bf a},{\bf b},{\bf p}\}} 
\big( {\bf y}_t \big) = {\rm const}_{{\bf a},{\bf p}} \Phi(t) , \quad
\psi^{(j)}_{{\bf p},{\rm Co}\{ {\bf a},{\bf b},{\bf p}\}} \big( {\bf y}_t \big) = {\rm const}^{(j)}_{{\bf a},{\bf p}} \Psi(t) , \\
&\varphi_{{\bf p},{\rm Co}\{ {\bf a},{\bf b},{\bf p}\}} 
\big( {\bf z}_t \big) = {\rm const}_{{\bf b},{\bf p}} \Phi(t) , \quad
\psi^{(j)}_{{\bf p},{\rm Co}\{ {\bf a},{\bf b},{\bf p}\}} \big( {\bf z}_t \big) = {\rm const}^{(j)}_{{\bf b},{\bf p}} \Psi(t) 
\end{aligned}
\eqno(\theIpvii.\theipvii)$$ %(1.7)$$
while for the points 
${\bf x}_t := (1-t){\bf a}+t{\bf b}$ on the edge 
${\rm Co}\{ {\bf a},{\bf b}\}$ 
we simply have
\ADT
\newcounter{Ipviii} \setcounter{Ipviii}{\value{SEC}}
\newcounter{ipviii} \setcounter{ipviii}{\value{THM}}
$$\begin{aligned}
&\varphi_{{\bf p},{\rm Co}\{ {\bf a},{\bf b},{\bf p}\}} ({\bf x}_t) = 
\psi^{(j)}_{{\bf p},{\rm Co}\{ {\bf a},{\bf b},{\bf p}\}} ({\bf x}_t) = 0, \\   
&\nabla \varphi_{{\bf p},{\rm Co}\{ {\bf a},{\bf b},{\bf p}\}} ({\bf x}_t) = 
\nabla \psi^{(j)}_{{\bf p},{\rm Co}\{ {\bf a},{\bf b},{\bf p}\}} ({\bf x}_t) = 
{\bf 0} . 
\end{aligned}
\eqno(\theIpviii.\theipviii)$$ %(1.8)$$
\rm

\noindent
In the sequel we call $\Phi,\Psi$ the {\it shape functions} of the
spline procedure $({\cal T},F) \mapsto f_{{\cal T},F}$ satisfying Postulate B. 
Notice that the requirements $\Phi(0)=\Phi^\prime(0) =\Psi(0)=\Psi^\prime(0)$
follow automatically from the 
%vanishing 
order condition $(\theIpv.\theipv)$ %(1.5) 
on the edge ${\rm Co}\{{\bf a,b}\}$.

At first glance, shape uniformity may seem an artificial requirement.
However,  
for a procedure satisfying Postulate A,
the geometrically natural property of  
being invariant with respect to homothetic transformations
$($maps $\RR \!\leftrightarrow\! \RR$ of the form 
${\bf x} \!\mapsto\! \mu {\bf xS} \!+\! {\bf w}$
with some orthogonal matrix ${\bf S} )$
implies Postulate B trivially. In our context we understand
invariance as follows:
\DEL{
\medskip
\noindent 
\bf Definition 1.9. \rm   
Let ${\bf G}:\RR^2\leftrightarrow\RR^2$, 
${\bf G(x)}= {\bf xA}+{\bf w}$ with some invertible matrix ${\bf A}$
be a surjective affine transformation of the plain. 
A spline procedure $\SS : ({\cal T},F) \mapsto f_{{\cal T},F}$ 
is {\it ${\bf G}$-invariant} if\footnote{$^*$}{Heuristics: 
${\bf G}$ should transfer spline functions 
constructed with the gradient data of any smooth function 
$h$ on ${\rm Vert}({\cal T})$ from ${\rm Dom}({\cal T})$ 
to the analogous objects with $h\circ{\bf G}^{-1}$ on 
${\rm Dom}({\bf G}({\cal T})) = {\bf G}({\rm Dom}({\cal T}))$.}, 
in any case we have
}%ENDDEL
given a surjective affine transformation 
${\bf G(x)}= {\bf xA}+{\bf w}$
with some invertible $2\!\times\!2$-matrix ${\bf A}$ 
of the plain,
the spline procedure $\SS : ({\cal T},F) \mapsto f_{{\cal T},F}$ 
is {\it ${\bf G}$-invariant} if
it transfers spline functions 
constructed with the gradient data of any smooth function 
$h$ on ${\rm Vert}({\cal T})$ from ${\rm Dom}({\cal T})$ 
to the analogous objects with $h\!\circ\!{\bf G}^{-1}$ on 
${\rm Dom}({\bf G}({\cal T})) \!=\! {\bf G}({\rm Dom}({\cal T}))$,
that is
\ADT
\newcounter{Ipx} \setcounter{Ipx}{\value{SEC}}
\newcounter{ipx} \setcounter{ipx}{\value{THM}} 
$$\begin{aligned}
&f_{{\cal T},F} \circ {\bf G}^{-1} =
f_{{\bf G}({\cal T}),{\bf G}^\sharp(F)} \\ 
&\qquad {\rm with} \quad
{\bf G}^\sharp \big( {\bf x},\chi,[\mu_1,\mu_2]\big) :=
\big( {\bf G}({\bf x}),\chi,[\mu_1,\mu_2]{\bf A}^{-1}\big) .
\end{aligned} 
\eqno(\NUM)$$ %(1.10)$$
As we shall see (Lemma 4.1), %\theIVpi.\theivpi
if Postulate A holds, we can formulate ${\bf G}$-invariance
in terms of the basic functions as follows:
\ADT
\newcounter{Ipxi} \setcounter{Ipxi}{\value{SEC}}
\newcounter{ipxi} \setcounter{ipxi}{\value{THM}}
$$\begin{aligned}
&\varphi_{{\bf G}({\bf p}),{\bf G}({\bf T})} \!=\! 
\varphi_{{\bf p,T}}\!\circ\! {\bf G}^{-1} ,\\
&\big[ \psi^{(1)}_{\bf G(p),G(T)} , \psi^{(2)}_{\bf G(p),G(T)}\big]
\!=\!
\big[ \psi^{(1)}_{\bf p,T}\!\circ\! {\bf G}^{-1} ,  
\psi^{(2)}_{\bf p,T}\!\circ\! {\bf G}^{-1} \big] {\bf A} .
\end{aligned}
\eqno(\NUM)$$ %(1.11)$$
It is worth to notice (Corollary 4.4) that 
$(\theIpxi.\theipxi)$ %$(1.11)$ 
cannot hold simultaneously for all
invertible matrices ${\bf A}$ and ${\bf w}\in \RR^2$. 
Thus there is no local linear spline procedure which is invariant 
under all invertible affine transformations and producing always 
${\cal C}^1$-smooth functions (i.e. functions being 
continuously differentiable also over the edges of mesh triangles)
functions.   

\medskip

Our aim in this paper is a parametric classification of the procedures
satisfying Postulates A,B, resulting in ${\cal C}^1$-smooth functions.
In particular we enumerate 
all the homothetically invariant linear local polynomial 
$\mathcal{C}^1$-spline interpolation procedues
from gradient data over triangular meshes.   
It is remarkable that
there is a unique one among them with lawest degree (degree 5) 
which turns out to be homothetically invariant.
From the view point of applications,
the results provide the complete list of hermition $\mathcal{C}^1$-splines
with shape uniformity over edges
from which one can choose the best fit one with respect to
various aspects.  
It is worth to relate the latter fact to a celebrated 
alternative local linear polynomial spline interpolation
procedure on the basis of
Zl\'amal-\v Zeni\v sek 
2-nd order triangular spline equations [\Zlam].
This relies upon the fact that,
given a triangular mesh with gradient and Hessian data at the vertices
and normal derivative values at edge middle points,
there is a unique fitting spline with 5th degree polynomials.
The 21 polynomial coefficients over any mesh triangle can be obtained
as the unique solution of a system of 21 straightforward linear equations
whose explicit formula was published recently [\Serg].
Though not stated in the sources, easily seen 
this kind of procedure has some homothetical invariance properties.
Hence 
it seems that our first order approach
with the shape conditions of Postulate B
provides a geometrically motivated
alternative to several problems discussed in [\Serg].
As mentioned earler and remarked also e.g. in [\Barn],
first order approches with a few (actually 9 in [\Barn]) free parameters
may have practical advantages versus higher oreder methods
due to the fact that data sampling can rarely support 
e.g. Hessian data (or even adequate guesses for them). 

Our arguments are based on the use of baricentic coordinates
associated with triangles instead of the usual Cartesian ones.
Applying Remark 3.2                 %\theIIIpii.\theiiipii\   %3.2 
to the difference of the first order
solution given in Theorem 2.3       %\theIIpii.\theiipii,  %%2.2, 
a way is opened to  
develop a new geometric approach to the system of
Zl\'amal-\v Zeni\v sek equations and its alternative variants 
which may have 
further independent theoretical and educational interest.

\bigskip\bigskip
\noindent
\ADS
{\bfs \theSEC. Main results}

\medskip
Recall that given a 
non-degenerate triangle ${\bf T} \subset \RR^2$ with
$\{ {\bf a,b,c} \} ={\rm Vert}({\bf T})$, 
the {\it normalized baricentric coordinates} 
of a point ${\bf x}$ 
are the terms of the
necessarily unique triple
$\big[ \lambda_{{\bf T}}^{{\bf a}}({\bf x}),
\lambda_{{\bf T}}^{{\bf b}}({\bf x}),
\lambda_{{\bf T}}^{{\bf c}}({\bf x}) \big] \in\RR^3$
such that
$${\bf x} = \lambda_{{\bf T}}^{{\bf a}}({\bf x}) {\bf a} +
\lambda_{{\bf T}}^{{\bf b}}({\bf x}) {\bf b} +
\lambda_{{\bf T}}^{{\bf c}}({\bf x}) {\bf c}, 
\quad
\lambda_{{\bf T}}^{{\bf a}}({\bf x}) +
\lambda_{{\bf T}}^{{\bf b}}({\bf x}) +
\lambda_{{\bf T}}^{{\bf c}}({\bf x}) =1 .
%\eqno(2.1)
$$
We reserve the symbols \
$\lambda_{{\bf T}}^{{\bf p}}$ \ 
as standard notation. 
It is well-known from elementary analytic plain geomertry 
[\Cox] that
$$\lambda_{\bf T}^{\bf p}({\bf x}) =
{ {\rm area}({\rm Co}\{ {\bf a,b,x}\})/
{\rm area}({\bf T}) } \qquad \big( {\bf x}\in{\bf T} \big)$$
thus normalized baricentric coordinates can easily be calculated 
by means of determinants or inner products with a $(\pi/2)$-rotation:
\ADT
\newcounter{IIpi} \setcounter{IIpi}{\value{SEC}}
\newcounter{iipi} \setcounter{iipi}{\value{THM}}
$$%\begin{aligned}
\lambda_{\bf T}^{\bf p}({\bf x}) \!=\! 
\frac{{\rm det}({\bf x \!-\! a},
{\bf x \!-\! b})}{{\rm det}({\bf p \!-\! a},{\bf p \!-\! b})} 
\frac{\big\langle {\bf (b \!-\! a)R} \big\vert {\bf x\!-\! a} \big\rangle}{
\big\langle {\bf (b-a)R} \big\vert {\bf p-a} \big\rangle} 
\ \ {\rm where} \ \ 
{\bf R} \!:=\!\left[\begin{matrix} 0 &1 \cr -1 &0 \end{matrix}\right] .
%\end{aligned} 
\eqno (\NUM)$$ %(2.1)$$ 
For later use we also introduce the abbreviating notations
\ADT
\newcounter{IIpiw} \setcounter{IIpiw}{\value{SEC}}
\newcounter{iipiw} \setcounter{iipiw}{\value{THM}}
$$x^{[j]}_{\bf p} \! := \! x^{[j]} \!-\! x^{[j]}({\bf p}), \
\xi^{\bf v}_{\bf p,a} \!:=\! \frac{
\langle {\bf v} \!-\! {\bf a \vert p\!-\! a} \rangle}{
\Vert {\bf p-a }\Vert^2 } , \
\overline{\xi}^{\bf v}_{\bf p,a} \!:=\!
\frac{\langle {\bf v\!-\! a} \vert {\bf (p\!-\! a)R} \rangle}{
\Vert {\bf p-a} \Vert^2 } .
\eqno(\NUM)$$  %(2.1^\prime)$$
As for geometric interpretation, 
$\xi^{\bf v}_{\bf p,a}$ resp. $\overline{\xi}^{\bf v}_{\bf p,a}$
are the affine coordinates of the point ${\bf v}$ with respect to the 
orthogonal frame $\big[ {\bf a, p, a \!+\!(p\!-\! a)R} \big]$
with origin ${\bf a}$ so that
${\bf v}= {\bf a} + \xi^{\bf v}_{\bf p,a} {\bf (p\!-\! a)} +
\overline{\xi}^{\bf v}_{\bf p,a} {\bf (p\!-\! a)R}$.
 
\bigskip
\noindent
\ADT\newcounter{IIpii} \setcounter{IIpii}{\value{SEC}}
\newcounter{iipii} \setcounter{iipii}{\value{THM}}
\bf Theorem \NUM.  %2.2. 
\it There is a unique local linear 
polynomial ${\cal C}^1$-spline 
procedure acting on triagular meshes with the property of 
uniform shape on vertices\footnote{\rm That is
satisfying Postulates A,B with 
$f_{{\cal T},F} \in {\cal C}^1\big({\rm Dom}({\cal T})\big)$.}
and having shape functions 
with minimal computational complexity.
Its shape functions are
%%\ADT
%%\newcounter{IIpiii} \setcounter{IIpiii}{\value{SEC}}
%%\newcounter{iipiii} \setcounter{iipiii}{\value{THM}}
$$\strut^*\!\Phi(t) = t^3 ( 10 - 15t + 6t^2 ), \qquad
\strut^*\!\Psi(t) = t^3 (t - 1) (4 - 3t ) .$$
%%\eqno(\NUM)$$  %(2.3)$$
The corresponding basic functions \  
$(\,$for a non-degenerate triangle \  
${\bf T} = {\rm Co}\{ {\bf a,b,p} \}$ 
with distinguished vertex ${\bf p} \,)$ 
have the form 
%%\ADT
%%\newcounter{IIpv} \setcounter{IIpv}{\value{SEC}}
%%\newcounter{iipv} \setcounter{iipv}{\value{THM}}
$$\hskip-2mm
\begin{aligned}
&\phantom{.}^*\!\varphi_{{\bf p},{\bf T}} \!=\! 
\strut^*\!\Phi\big( \lambda_{\bf T}^{\bf p} \big) + 
30 \,[\lambda_{\bf T}^{\bf p}]^2 
\lambda_{\bf T}^{\bf a} \lambda_{\bf T}^{\bf b} 
\Big[ 
\xi^{\bf b}_{\bf p,a}
\lambda_{\bf T}^{\bf b} +
\xi^{\bf a}_{\bf p,b}
\lambda_{\bf T}^{\bf a} \Big] ,\\
&\phantom{.}^*\!\psi^{(j)}_{{\bf p}, {\bf T}} \!=\!
\frac{\strut^*\!\Psi\big(\lambda_{\bf T}^{\bf p}\big) }{
\lambda_{\bf T}^{\bf p} - 1} 
x^{[j]}_{\bf p} + 12 [\lambda_{\bf T}^{\bf p}]^2 
\lambda_{\bf T}^{\bf a} \lambda_{\bf T}^{\bf b} \left[ 
\xi^{\bf b}_{\bf p,a}
x^{[j]}_{\bf p}({\bf b}) \lambda_{\bf T}^{\bf b} +
\xi^{\bf a}_{\bf p,b}
x^{[j]}_{\bf p}({\bf a}) \lambda_{\bf T}^{\bf a} \right] \! . 
\end{aligned}$$
%%\eqno(\NUM)$$  %(2.5)$$

\ADT
\newcounter{IIpvi} \setcounter{IIpvi}{\value{SEC}}
\newcounter{iipvi} \setcounter{iipvi}{\value{THM}}
\noindent
\bf Theorem \NUM.  %2.6. 
\it  A spline procedure acting on triangular meshes 
and satisfying Postulates A,B produces ${\cal C}^1$-smooth splines
if and only if
its shape functions are of the form
$$\Phi(t) = \strut^*\!\Phi(t) +t^3(1-t)^3 \Phi_1(t), \qquad
\Psi(t) = \strut^*\!\Psi(t) + t^3(1-t)^3 \Psi_1(t) 
\eqno(2.7)$$ 
and the basic functions
$($for a non-degenerate triangle
${\bf T}={\rm Co}\{ {\bf a,b,p}\}$ with distinguished vertex ${\bf p} )$
can be written
in terms of the modified shape function
%%\ADT
%%\newcounter{IIpvw} \setcounter{IIpvw}{\value{SEC}}
%%\newcounter{iipvw} \setcounter{iipvw}{\value{THM}}
$$\Theta(t) \!:=\! \Psi(t)/(t \!-\! 1)$$
%%\eqno(\NUM)$$  %(2.5^\prime)$$ 
and the rotation matrix ${\bf R}$ in $(\theIIpi.\theiipi)$ %$(2.1)$ 
as
\ADT
\newcounter{IIpviii} \setcounter{IIpviii}{\value{SEC}}
\newcounter{iipviii} \setcounter{iipviii}{\value{THM}}
$$\begin{aligned}
&\varphi_{\bf p, T} =
\Phi\big(\lambda_{\bf T}^{\bf p}\big) + [\lambda_{\bf T}^{\bf p}]^2 
\lambda_{\bf T}^{\bf a} \lambda_{\bf T}^{\bf b} 
P^{{\bf p}}_{\bf a,b}\big(
\lambda_{\bf T}^{\bf b},\lambda_{\bf T}^{\bf a}\big), \\
&\psi^{(j)}_{\bf p,T} =
\Theta\big(\lambda_{\bf T}^{\bf p}\big)
x_{{\bf p}}^{[j]} + [\lambda_{\bf T}^{\bf p}]^2 
\lambda_{\bf T}^{\bf a} \lambda_{\bf T}^{\bf b} 
Q^{j,{\bf p}}_{\bf a,b} \big(
\lambda_{\bf T}^{\bf b},\lambda_{\bf T}^{\bf a}\big)
\end{aligned}
\eqno(\NUM)$$  %(2.8)$$
where
\ADT
\newcounter{IIpix} \setcounter{IIpix}{\value{SEC}}
\newcounter{iipix} \setcounter{iipix}{\value{THM}}
$$\begin{aligned}
&P^{\bf p}_{\bf a,b}(s,t) = 
s \Big\{ 
\xi^{\bf a}_{\bf p,b}
\frac{\Phi^\prime(1-s) }{ (1-s)^2s^2} +
\overline{\xi}^{\bf a}_{\bf p,b}
k_{\bf b}^{0,{\bf p}}(s) \Big\} +\\ 
&\hskip25mm 
+ t \Big\{ 
\xi^{\bf b}_{\bf p,a}
\frac{\Phi^\prime(1-t) }{ (1-t)^2t^2}
\overline{\xi}^{\bf b}_{\bf p,a}
k_{\bf a}^{0,{\bf p}}(t) \Big\}  
+st R^{0,{\bf p}}_{\bf a,b}(s,t) , \\
&Q^{j,{\bf p}}_{\bf a,b} (s,t) = 
s \Big\{ 
\xi^{\bf a}_{\bf p,b} 
x^{[j]}_{\bf p}\!({\bf b})
\frac{\Theta^\prime(1 \!-\! s) }{ s(1\!-\! s)^2} +
\overline{\xi}^{\bf a}_{\bf p,b}  
k_{\bf b}^{j,{\bf p}}(s) \Big\} + \\ 
& \hskip25mm + 
t \Big\{ 
\xi^{\bf b}_{\bf p,a}  
x^{[j]}_{\bf p}\!({\bf a})
\frac{\Theta^\prime(1 \!-\!t) }{ t(1 \!-\! t)^2} +
\overline{\xi}^{\bf b}_{\bf p,a}
k_{\bf a}^{j,{\bf p}}(t) \Big\} + 
st R^{j,{\bf p}}_{\bf a,b}(s,t) 
\end{aligned}
\eqno(\NUM)$$  %(2.9)$$
with the following free options in  $(\theIIpviii.\theiipviii)$ resp. $(\theIIpix,\theiipix)$:  
%$(2.7)$ resp. $(2.9)$:
\medskip
\begin{itemize}
\item[{\rm (i)}]
$\Phi_1,\Psi_1:[0,1]\to\RR$ are arbitrary polynomial functions,

\item[{\rm (ii)}]
$({\bf p,q}) \mapsto k_{\bf q}^{i,{\bf p}}$ $(i=0,1,2)$
are arbitrary maps assigning polynomial functions $\RR\to\RR$
to pairs of distinct points,

\item[{\rm (iii)}]
$({\bf p,q,r}) \mapsto 
R^{i,{\bf p}}_{\bf q,r}$ $(i=0,1,2)$ 
are arbitrary maps assigning polynomial functions  $\RR^2\to\RR$
to triples of non-collinear points
with the symmetry
$ R^{i,{\bf p}}_{\bf q,r}(s,t) \equiv R^{i,{\bf p}}_{\bf r,q}(t,s)$.
\end{itemize}

\ADT
\newcounter{IIpx} \setcounter{IIpx}{\value{SEC}}
\newcounter{iipx} \setcounter{iipx}{\value{THM}}
\bigskip
\noindent
\bf Remark \NUM.  %2.10. 
\rm (i) Actually, Theorem \theIIpii.\theiipii\  %2.2 
is simply a corollary
of Theorem \theIIpvi.\theiipvi\ %2.6 
by setting the options ${\rm (i)}\!-\!{\rm (iv)}$
to $0$. We emphasize it for its potential practical and educational use.

\smallskip

(ii)  
The formally rational expressions 
in $(\theIIpviii.\theiipviii \!-\! \theIIpix.\theiipix)$  %$(2.7)$ 
are polynomials.
Indeed, $\Phi^\prime(1-t) / \big[ t^2(1-t)^2 \big] = 
30 - 3(1-2t)\Phi_1(1-t) + t(1-t) \Phi_1^\prime(1-t)$, 
resp. \  
$\Psi(t)/(t-1) = t^3[(4-3t) - (1-t)^2\Psi_1(t)] , \ \ 
\Theta^\prime(1-t)/ [t(1-t)^2] = 
12 + (2-5t)\Psi_1(1-t) - t(1-t)\Psi_1^\prime(1-t)$.  

\smallskip

(iii) 
$\lambda_{\bf T}^{\bf p},\lambda_{\bf T}^{\bf a},\lambda_{\bf T}^{\bf b}$
are the affine functions  
determined
by the properties
${\rm Line}\{\! {\bf a,\! b}\!\} \!\!=\!\! 
\big( \lambda_{\bf T}^{\bf p} \!\!=\!\!0 \big)$,
${\rm Line}\{\! {\bf b,\! p} \!\} \!\!=\!\! 
\big( \lambda_{\bf T}^{\bf a} \!\!=\!\! 0 \big)$,  
${\rm Line}\{\! {\bf a,\! p} \!\} \!\!=\!\! 
\big( \lambda_{\bf T}^{\bf b} \!\!=\!\!0 \big)$, 
$\phantom{\strut^{\int}_{\int}}$\hskip-3mm  
$\lambda_{\bf T}^{\bf p}({\bf p}) \!\!=\!\! 
\lambda_{\bf T}^{\bf a}({\bf a}) \!\!=\!\!
\lambda_{\bf T}^{\bf b}({\bf b}) \!\!=\!\! 1$.
For the parametrized edge points in $(\theIpvii.\theipvii)$ %$(1.7)$ 
we have 
\begin{figure}[h]
\begin{center}
%\hfill\break
\vtop{
\hbox{$\phantom{\strut^{\int}}$}
\hbox{$\phantom{\strut^{\int}}$ 
$\lambda_{\bf T}^{\bf p}({\bf x}_t) =
\lambda_{\bf T}^{\bf a}({\bf z}_t) = 
\lambda_{\bf T}^{\bf b}({\bf y}_t) \equiv 0$,}
\hbox{$\phantom{\strut^{\int}}$ 
$\lambda_{\bf T}^{\bf p}({\bf y}_t) =
\lambda_{\bf T}^{\bf p}({\bf z}_t) =
\lambda_{\bf T}^{\bf b}({\bf x}_t) \equiv t$,}
\hbox{$\phantom{\strut^{\int}_{\int}}$
$\lambda_{\bf T}^{\bf a}({\bf x}_t) = 
\lambda_{\bf T}^{\bf a}({\bf y}_t) =
\lambda_{\bf T}^{\bf b}({\bf z}_t) \equiv 1-t$.}
\hbox{$\phantom{\strut^{\int}}$}}
\qquad
\vskip-2.7cm
\begin{tikzpicture}[x=5mm, y=5mm]; % rotate=30]
\node at (-8,0) {};
\draw (5,0) -- (10,0) -- (7,3) -- (5,0);
\node at (5,0) {{\hbox{$\bullet$}}};
\node at (10,0) {{\hbox{$\bullet$}}};
\node at (7,3) {{\hbox{$\bullet$}}};
\node at (4.5,0.2) {{\bf a}};
\node at (10.4,0.2) {{\bf b}};
\node at (9.5,2.9) {\hbox{${\bf p=y_1=z_1}$}};
\node at (6.5,0) {{\hbox{$\bullet$}}};
\node at (6.5,0.4) {{\hbox{${\bf x}_t$}}};
\node at (5.6,0.9) {{\hbox{$\bullet$}}};
\node at (5.0,0.9) {{\hbox{${\bf y}_t$}}};
\node at (9.1,0.9) {{\hbox{$\bullet$}}};
\node at (9.7,0.9) {{\hbox{${\bf z}_t$}}};
\end{tikzpicture}
\end{center}
\end{figure}
%\break 
On the other hand 
$x^{[j]}({\bf y}_t) = 
(1-t) x^{[j]} ({\bf a-p}) = (1-t) x^{[j]}_{\bf p}({\bf a})$
resp.
$x^{[j]}({\bf z}_t) = 
(1-t) x^{[j]} ({\bf b-p}) = (1-t) x^{[j]}_{\bf p}({\bf b})$.  
Hence, with the formulas $(\theIIpviii.\theiipviii)$,  %$(2.8)$, 
the shape conditions $(\theIpvii.\theipvii)$ %(1.7) 
hold automatically
with 
${\rm const}_{{\bf a},{\bf p}} = {\rm const}_{{\bf b},{\bf p}} =1$ and
${\rm const}^{(j)}_{{\bf a},{\bf p}} = x^{[j]}({\bf p-a})$ 
resp.  
${\rm const}^{(j)}_{{\bf b},{\bf p}} = x^{[j]}({\bf p-b})$,
furthermore also $(\theIpviii.\theipviii)$ %(1.8) 
is fulfilled.

\smallskip

(iv)
One can check with symbolic computer algebra that
all the spline procedures described in Theorem \theIIpvi.\theiipvi\   %2.6 
produce
${\cal C}^1$-functions. It suffices to establish only
that, given any two adjacent non-degenerate triangles 
${\bf T}:={\rm Co}\{ {\bf p,a,b}\}$ resp.
$\widetilde{\bf T}:={\rm Co}\{ {\bf p,a},\widetilde{\bf p}\}$
with common edge ${\rm Co}\{ {\bf p,a}\}$ and distinguished point
${\bf p}$, the gradient vectors of the basic functions
$\varphi_{\bf p,T}, \psi^{(j)}_{\bf p,T}$ coincide with those of
$\varphi_{{\bf p},\widetilde{\bf T}}, 
\psi^{(j)}_{{\bf p},\widetilde{T}}$ at the points
${\bf y}_t =(1-t){\bf a}+t{\bf b}$.
Indeed, hence it follows that  
the unit spline functions
$f_{{\cal T},F_{\bf p}^i}$ 
$\big( {\bf p}\!\in\!{\rm Vert}({\cal T}),\, i\!=\!0,1,2\big)$
corresponding to the gradient data
\DEL{
$F^i_{\bf p}:= 
\big\{ [{\bf p},1,{\bf 0}]$ if $i\!=\! 0$, 
$[{\bf p},0,{\bf e}]^{[i]}$ if $i\!=\!1,2 \big\} \cup
\big\{ [{\bf q},0,{\bf 0}]: 
{\bf p}\ne{\bf q}\in {\rm Vert}({\cal T}) \big\}$
}%ENDDEL
$F^0_{\bf p}:= 
\big\{ [{\bf p},1,{\bf 0}],
[{\bf q},0,{\bf 0}]: 
{\bf q}\in {\rm Vert}({\cal T})\setminus\{ {\bf p}\} \big\}$
resp.
$F^j_{\bf p}:= \big\{ [{\bf p},0,{\bf e}^{[j]}],
[{\bf q},0,{\bf 0}]: 
{\bf q}\in {\rm Vert}({\cal T})\setminus\{ {\bf p}\} \big\}$
$(j\!=\! 1,2)$
are continuously differentiable.
%%A related MAPLE code can be found in [\St]. 

\ADT
\newcounter{IIpxi} \setcounter{IIpxi}{\value{SEC}}
\newcounter{iipxi} \setcounter{iipxi}{\value{THM}}
\bigskip
\noindent
\bf Theorem \NUM.  %2.11.  
\it A ${\cal C}^1$-spline procedure 
$\hbox{\Goth S}$ described in Theorem \theIIpvi.\theiipvi\  %$2.6$
in the form $(2.7\!\!-\!\! 9)$ is isometry-invariant
if and only if 
$k^{i,{\bf p}}_{\bf c}(t) \equiv 0$
for all  
$i \!=\!0,1,2$ and ${\bf p\!\ne\! c}\!\in\!\RR^2$
furthermore the higher terms $R^{i,{\bf p}}_{\bf a,b}$ in $(2.9)$
transform as
$R^{0,{\bf G\!(p)}}_{\bf G\!(a),G\!(b)} = 
R^{0,{\bf p}}_{\bf a,b} \phantom{\strut^{\strut}_{\strut}}$ resp.
$\big [R^{1,{\bf G\!(p)}}_{\bf G\!(a),G\!(b)}, 
R^{2,{\bf G\!(p)}}_{\bf G\!(a),G\!(b)} \big] = 
\big[ R^{1,{\bf p}}_{\bf a,b} , R^{2,{\bf p}}_{\bf a,b} \big] {\bf A}$
whenever ${\bf G}: {\bf x}\mapsto {\bf w + x A}$ is an isometry.

\bigskip
\bigskip
\noindent
\ADS
{\bfs 3. Proof of Theorem \theIIpvi.\theiipvi}  %2.6}

\medskip
%\noindent
\rm
Henceforth we consider an arbitrarily fixed 
procedure
$\SS : ({\cal T},F) \!\mapsto\! f_{{\cal T},F}$
which satisfies Postulates A,B 
and produces continuous but not necessarily
continuously differentiable functions. 
We reserve the notations 
$\varphi_{{\bf p},{\bf T}}$, $\psi^{(j)}_{{\bf p},{\bf T}}$
resp. $\Phi,\Psi$ 
for the basic functions resp. shape functions
as established in Section \theSi. %Section 1.  
In accordance with $(\theIpvi.\theipvi)$ %(1.6),
we can write
\ADT
\newcounter{IIIpi} \setcounter{IIIpi}{\value{SEC}}
\newcounter{iiipi} \setcounter{iiipi}{\value{THM}} 
$$\Phi(t) = t^2(3-2t) + t^2(1-t)^2 \Phi_0(t), \qquad
\Psi(t) = t^2(t-1) + t^2(1-t)^2\Psi_0(t)
\eqno(\NUM)$$   %(3.1)$$
and \ $\Theta(t) = t^2 + t^2(t-1)\Psi_0(t)$ \ 
with suitable polynomials $\Phi_0,\Psi_0$.

Next we are going to express
the constraints $(\theIpv.\theipv)$, %(1.5),
$(\theIpvii.\theipvii\!-\!\theipviii)$ %(1.7-8) 
on the basic functions in terms of
$\Phi,\Psi$ and 
baricentric coordinates. 
To this aim, we recall the following folklore fact 
from elementary algebraic geometry
relating the root curves with a product decomposition
of multivariate polynomials 
which is an easy consequence of 
B\'ezout's Theorem [\Fulton,\Hamilton].

\ADT
\newcounter{IIIpii} \setcounter{IIIpii}{\value{SEC}}
\newcounter{iiipii} \setcounter{iiipii}{\value{THM}}
\bigskip
\noindent
{\bf Remark \NUM}(i) %3.2}(i)
\ If \ 
${\bf L}_0,{\bf L}_1,\ldots,{\bf L}_m$
are distinct straight lines such that 
${\bf L}_k = \big(\ell_k =0\big)$ with the
affine functions (i.e. polynomials of first degree)
$\ell_k: \RR^2\to\RR$ $(k=1,\ldots,m)$
then a polynomial $\RR^2\to\RR$ is divisable with
$\prod_{k=0}^m \ell_k^{\nu_k}$ 
if and only if, for any index $k$, it vanishes in order $\nu_k$
at the points of ${\bf L}_k$.
In particular, 
given a non-degenerate triangle
${\bf T}:={\rm Co}\{ {\bf a,b,p}\}$,
a polynomial $Q: \RR^2 \to \RR$ of two variables 
has the form
\ADT
%%\newcounter{IIIpiii} \setcounter{IIIpiii}{\value{SEC}}
%%\newcounter{iiipiii} \setcounter{iiipiii}{\value{THM}}
$$Q = 
[\lambda_{\bf T}^{\bf p}]^{\nu_0} [\lambda_{\bf T}^{\bf a}]^{\nu_1} 
[\lambda_{\bf T}^{\bf b}]^{\nu_2} $$
%%\eqno(\NUM)$$  %(3.3)$$
for some polynomial $q:\RR^2\to\RR$ if and only if
it vanishes in 
order $\nu_0$ at the points of ${\rm Line}\{{\bf a},{\bf b}\}$,
order $\nu_1$ at ${\rm Line}\{{\bf p},{\bf b}\}$
and
order $\nu_2$ at ${\rm Line}\{{\bf p},{\bf a}\}$,
respectively.\footnote{\rm $Q$ vanishes in order $\nu$ 
at the point $[x_0,y_0]$ if
$\frac{\partial^{k+m} }{ \partial x^k \partial y^m} Q(x_0,y_0) =0$
whenever $k+m<\nu$.}  

\medskip
(ii) \ If $Q:\RR^2\to\RR$ is a polynomial of two variables, we can write
%%\ADT
%%\newcounter{IIIpiiiw} \setcounter{IIIpiiiw}{\value{SEC}}
%%\newcounter{iiipiiiw} \setcounter{iiipiiiw}{\value{THM}}
$$\begin{aligned}
Q&(x,y) = Q(0,0) + x q_1(x) + y q_2(y) + xy q_3(x,y) \quad {\rm where}\\ 
&q_1(x) := [Q(x,0)]-Q(0,0)/x, \quad q_2(y):=[Q(0,y)-Q(0,0)]/y,\\  
&q_3(x,y) := \big[ Q(x,y) - [Q(0,0)+xq_1(x)+yq_2(y)]\big]/(xy) 
\end{aligned}$$
%%\eqno(\NUM)$$  %(3.3^\prime)$$
are well-defined polynomials in one resp. two variables.
We shall call the $\RR^2$-polynomial \ 
$Q_0(x,y):=Q(0,0) + x q_1(x) + y q_2(y)$
of first degree
the {\it principal part} of $Q$.
 
\ADT
\newcounter{IIIpiv} \setcounter{IIIpiv}{\value{SEC}}
\newcounter{iiipiv} \setcounter{iiipiv}{\value{THM}}
\bigskip
\noindent
{\bf Lemma \NUM.}  %3.4.} 
\it The basic functions \ 
$\varphi_{\bf p,T}, 
\psi^{(j)}_{\bf p,T}$ \ 
for \ ${\bf T}={\rm Co}\{ {\bf a,b,p}\}$ \ 
have the form
%%\ADT
%%\newcounter{IIIpv} \setcounter{IIIpv}{\value{SEC}}
%%\newcounter{iiipv} \setcounter{iiipv}{\value{THM}} 
$$\begin{aligned} 
&\varphi_{\bf p,T} = \Phi (\lambda_{\bf T}^{\bf p}) + 
[\lambda_{\bf T}^{\bf p}]^2 \lambda_{\bf T}^{\bf a} \lambda_{\bf T}^{\bf b} \,
{\rm Pol}(\lambda_{\bf T}^{\bf b},\lambda_{\bf T}^{\bf a}) , \\
&\psi^{(j)}_{\bf p,T} = \Theta(\lambda_{\bf T}^{\bf p}) x_{{\bf p}}^{[j]} + 
[\lambda_{\bf T}^{\bf p}]^2 \lambda_{\bf T}^{\bf a} \lambda_{\bf T}^{\bf b} \,
{\rm Pol}(\lambda_{\bf T}^{\bf b},\lambda_{\bf T}^{\bf a})
%\quad {\rm where} \quad \Theta(t):= \frac{\Psi(t) }{ t-1}
\end{aligned}$$
%%\eqno(\NUM)$$  %(3.5)$$
in terms of the 
baricentric coordinates $(\theIIpi.\theiipi)$,  %$(2.1)$, 
the shape functions $\theIIIpi.\theiiipi$,  %$(3.1)$,
$\Theta(t) \!:=\!\Psi(t)/(t\!-\! 1)$
and 
with suitable polynomials 
of two variables.

\bigskip
\noindent
\bf Proof. \rm
Fix any triangle ${\bf T}:={\rm Co}\{ {\bf a,b,p}\}$.
As mentioned, necessarily  $(\theIIIpi.iiipi)$  %(3.1) 
holds and $\Theta$ is a polynomial.
Consider the functions
$$f := \Phi(\lambda_{\bf T}^{\bf p}), \qquad 
g^{(j)} := \Theta(\lambda_{\bf T}^{\bf p}) \cdot x^{[j]}_{{\bf p}} .$$ 
Along the edge ${\rm Co}\{{\bf a,p}\}$, at the points \
${\bf y}_t := (1-t){\bf a} + t{\bf p}$ \ 
we have \ 
$\lambda_{\bf T}^{\bf p}({\bf y}_t) =t$, 
$\lambda_{\bf T}^{\bf b}({\bf y}_t) =0$, 
$\lambda_{\bf T}^{\bf a}({\bf y}_t) =
[1 -\lambda_{\bf T}^{\bf a} ({\bf y}_t) - \lambda_{\bf T}^{\bf b}({\bf y}_t)
=1-t$.
Observe that the functions $f,g^{(j)}$ suit the 
shape uniformity conditions because
$$\begin{aligned}
f({\bf y}_t) \!=\! \Phi(t), \quad 
g^{(j)} ({\bf y}_t) \!&=\! \Theta(t) 
\big\langle {\bf e}^{[j]} \big\vert {\bf y}_t -{\bf p} \big\rangle \!=\!
\Theta(t) (1-t) \big\langle {\bf e}^{[j]} \big\vert {\bf a-p} \big\rangle \!=\\
&= \big\langle {\bf e}^{[j]} \big\vert {\bf p-a} \big\rangle \Psi(t)
\end{aligned}$$
and since $f,g^{(j)}$ are polynomial multiples of 
$[\lambda_{\bf T}^{\bf p}]^2$.
Also, since ${\bf y}_1={\bf p}$, 
$f({\bf p}) =\Phi(\lambda_{\bf T}^{\bf p}({\bf y}_1))=\Phi(1)=1$ \  
and
$$\begin{aligned}
\nabla f({\bf p}) &= 
\Phi^\prime\big(\lambda_{\bf T}^{\bf p}({\bf y}_1)\big) 
\nabla \lambda_{\bf T}^{\bf p}({\bf y}_1) =
0\cdot \nabla \lambda_{\bf T}^{\bf p}({\bf y}_1) ={\bf 0} ,\\
\nabla g^{(j)} ({\bf p}) &= \nabla_{{\bf x}={\bf y}_1}
\Big[ \Theta\big(\lambda_{\bf T}^{\bf p}({\bf x})\big)x^{[j]}_{{\bf p}}({\bf x}) \Big] =\\
&= x^{[j]}_{{\bf p}}({\bf p}) \nabla_{{\bf x}={\bf y}_1} 
\Theta\big(\lambda_{\bf T}^{\bf p}({\bf x})\big) + 
\Theta\big( \lambda_{\bf T}^{\bf p}({\bf p})\big) \nabla_{{\bf x}={\bf y}_1} x^{[j]}({\bf p}) =\\
&= 0 \cdot \Theta^\prime(1) \nabla \lambda_{\bf T}^{\bf p}({\bf y}_1) +
\Theta(1) {\bf e}^{[j]} =
{\bf e}^{[j]} .
\end{aligned}$$
Therefore \ the difference functions \  
$\varphi_{\bf pT} \!-\! f$ \ and \ $\psi^{(j)}_{\bf p,T} \!-\! g^{(j)}$ 
vanish on the edge \ ${\rm Co}\{{\bf a,p}\}$ \ 
of the triangle ${\bf T}$.
Similar arguments with the points ${\bf z}_t:=(1-t){\bf b}+t{\bf p}$ 
show that 
$\varphi_{\bf p,T} \!-\! g$ and $\psi^{(j)}_{\bf p,T} \!-\! g^{(j)}$ 
vanish on ${\rm Co}\{{\bf b,p}\}$.
By $(\theIpviii.\theipviii)$  %$(1.8)$, 
their gradients also vanish on the edge 
${\rm Co}\{ {\bf a,b}\} = \big( \lambda_{\bf T}^{\bf p}=0\big)$. Hence
(cf. Remark  \theIIIpii.\theiiipii)  %3.2) 
they are
polynomial multiples of 
$[\lambda_{\bf T}^{\bf p}]^2\lambda_{\bf T}^{\bf a}\lambda_{\bf T}^{\bf b}$,
say 
$\varphi_{\bf pT} = f + 
[\lambda_{\bf T}^{\bf p}]^2\lambda_{\bf T}^{\bf a}\lambda_{\bf T}^{\bf b} 
\Pi^{(0)}_{\bf p,T}$ and 
$\psi^{(j)}_{\bf p,T} = g^{(j)} + 
[\lambda_{\bf T}^{\bf p}]^2\lambda_{\bf T}^{\bf a}\lambda_{\bf T}^{\bf b} 
\Pi^{(j)}_{\bf p,T}$, respectively.
Since $\lambda^{\bf a}_{\bf T},\lambda^{\bf b}_{\bf T}$ 
are linearly independent affine functionals,
the mapping 
$\Lambda^{\bf p}_{\bf a,b} :{\bf x}\mapsto \big[
\lambda^{\bf b}_{\bf T}({\bf x}),\lambda^{\bf a}_{\bf T}({\bf x}) \big]$
is an affine coordinatization on the plain $\RR^2$.
Thus we can express each term $\Pi^{(i)}_{\bf p,T}$ as a polynomial
of the coordinates $\Lambda^{\bf p}_{\bf a,b}$
which completes the proof.

\ADT
\newcounter{IIIpxvii} \setcounter{IIIpxvii}{\value{SEC}}
\newcounter{iiipxvii} \setcounter{iiipxvii}{\value{THM}}
\bigskip
\noindent 
\bf Notation \NUM.  %3.17. 
\rm  
For later convenience, without danger of confusion, 
we introduce the unifying context-free notations
$$\lambda^{\bf p}_{\bf a,b} := \lambda^{\bf p}_{{\rm Co}\{ {\bf a,b,p}\}},
\quad
f^{0,{\bf p}}_{\bf a,b} := \varphi_{{\bf p},{\rm Co}\{ {\bf a,b,p}\}}, 
\quad
f^{j,{\bf p}}_{\bf a,b} := 
\psi^{(j)}_{{\bf p},{\rm Co}\{ {\bf a,b,p}\}} \quad (j\!=\! 1,2) .$$
Furthermore, in view of Lemma \theIIIpiv.\theiiipiv,  %3.4, 
we shall write
\ADT
\newcounter{IIIpxviiw} \setcounter{IIIpxviiw}{\value{SEC}}
\newcounter{iiipxviiw} \setcounter{iiipxviiw}{\value{THM}} 
$$f_{{\bf a,b}}^{i,{\bf p}} =
\Phi^{[i]} (\lambda_{\bf a,b}^{\bf p}) x^{[i]}_{{\bf p}} + 
[\lambda_{\bf a,b}^{\bf p}]^2 \lambda_{\bf a,p}^{\bf b} 
\lambda_{\bf b,p}^{\bf a} 
P^{i,{\bf p}}_{\bf a,b}(\lambda_{\bf a,p}^{\bf b},\lambda_{\bf b,p}^{\bf a})
\quad (i=0,1,2)   
\eqno(\NUM)$$  %(3.17^\prime)$$
where
$$\Phi^{[0]} := \Phi, \quad \Phi^{[1]} := \Phi^{[2]} := \Theta,
\quad x^{[0]}_{\bf p} : {\bf x}\mapsto 1 \ \ 
{\rm with} \ \ {\bf e}^{[0]} := \nabla x^{[0]}_{\bf p} = {\bf 0}$$
and the terms  
$P^{i,{\bf p}}_{\bf a,b}$ $(i\!=\! 0,1,2)$ 
are polynomials with coefficients 
depending %analytically 
on the ordered tuple $(i,{\bf p,a,b})$.
Notice that necessarily
\ADT
\newcounter{IIIpxviiww} \setcounter{IIIpxviiww}{\value{SEC}}
\newcounter{iiipxviiww} \setcounter{iiipxviiww}{\value{THM}}
$$P^{i,{\bf p}}_{\bf a,b} (s,t) \equiv P^{i,{\bf p}}_{\bf b,a}(t,s)
\eqno(\NUM)$$  %(3.17^{\prime\prime})$$
due to the trivial index symmetries \ 
$\lambda_{\bf u,v}^{\bf w} \equiv \lambda_{\bf v,u}^{\bf w}$ \ 
and \ 
$f_{{\bf a,b}}^{i,{\bf p}} \equiv f_{{\bf b,a}}^{i,{\bf p}}$.

\ADT
\newcounter{IIIpxv} \setcounter{IIIpxv}{\value{SEC}}
\newcounter{iiipxv} \setcounter{iiipxv}{\value{THM}}
\bigskip
\noindent
{\bf Lemma \NUM.} \it  %3.15.} \it
We have $f_{{\cal T},F} \in {\cal C}^1\big({\rm Dom}({\cal T})\big)$
for every triangular mesh with arbitrary gradient data 
if and only if 
\ADT
\newcounter{IIIpxvw} \setcounter{IIIpxvw}{\value{SEC}}
\newcounter{iiipxvw} \setcounter{iiipxvw}{\value{THM}}
$${\bf b} \mapsto \nabla f_{{\bf a,b}}^{i, {\bf p}}({\bf y}) \equiv 
{\rm const}_{\bf p,a,y} \quad 
\hbox{\rm for fixed ${\bf p\ne a}$ and ${\bf y}\in{\rm Co}\{ {\bf p,a}\}$}. 
\eqno(NUM)$$   %(3.15^\prime)$$ 
\noindent
{\bf Proof.} \rm Given any non-degenerate triangle 
${\bf T}={\rm Co}\{ {\bf a,b,p}\}$,
By construction, for the points
${\bf x}_t := (1-t){\bf a} + t {\bf b}$, 
${\bf y}_t := (1-t){\bf a} + t {\bf p}$
and
${\bf z}_t := (1-t){\bf b} + t {\bf p}$
on the edges of the triangle ${\bf T}$ 
we have \ 
$f^{i,{\bf p}}_{\bf a,b}({\bf x}_t) = 0$
independently of ${\bf p}$, \ 
$f^{i,{\bf p}}_{\bf a,b}({\bf y}_t) = \Phi^{[i]}(t)$
independently of ${\bf b}$ and
$f^{i,{\bf p}}_{\bf a,b}({\bf z}_t) = \Phi^{[i]}(t)$
independently of ${\bf a}$.
Thus the shape conditions are automatic from $(\theIIIpxviiw.\theiiipxviiw)$.  %$(3.17^\prime)$.
Moreover, given any triangle $\widetilde{\bf T}$ 
with a common edge but disjoint interior to ${\bf T}$,  
the functions pairs 
$\varphi_{\bf p,T},\varphi_{\bf p,\widetilde{T}}$ 
resp.
$\psi^{(j)}_{\bf p,T},\psi^{(j)}_{\bf p,\widetilde{T}}$
touch continuosly.
The analogous necessary and sufficent condition for a ${\cal C}^1$-smooth
touching is that the gradient pairs
$\nabla\varphi_{\bf p,T},\nabla\varphi_{\bf p,\widetilde{T}}$ 
resp.
$\nabla\psi^{(j)}_{\bf p,T},\nabla\psi^{(j)}_{\bf p,\widetilde{T}}$
coincide on the common edge:
\ADT
\newcounter{IIIpxiii} \setcounter{IIIpxiii}{\value{SEC}}
\newcounter{iiipxiii} \setcounter{iiipxiii}{\value{THM}}
$$\begin{aligned}
{\rm (i)} \ \ &\nabla f^{i, {\bf p}}_{{\bf a,b}} ({\bf x}) =
\nabla f^{i, \widetilde{\bf p}}_{{\bf a,b}} ({\bf x}) \hskip10mm
{\rm if} \  {\bf x}\in{\rm Co}\{ {\bf a,b} \} \!=\!
{\bf T}\!\cap\! {\rm Co}\{ {\bf a,b},\widetilde{\bf p}\},
\\ 
{\rm (ii)} \ \ &\nabla f^{i, {\bf p}}_{{\bf a,b}} ({\bf y}) =
\nabla f^{i, \widetilde{\bf b}}_{{\bf a,p}} ({\bf y}) \hskip10mm
{\rm if} \ {\bf y}\in{\rm Co}\{ {\bf a,p} \}  \!=\!
{\bf T}\!\cap\!{\rm Co}\{ {\bf a},\widetilde{\bf b},{\bf p}\}, \\
{\rm (iii)} \ \ &\nabla f^{i, {\bf p}}_{{\bf a,b}} ({\bf z}) =
\nabla f^{i, \widetilde{\bf a}}_{{\bf b,p}} ({\bf z}) \hskip10mm
{\rm if} \ {\bf z}\in{\rm Co}\{ {\bf b,p} \} \!=\!
{\bf T}\!\cap\!{\rm Co}\{ \widetilde{\bf a},{\bf b,p}\} .
\end{aligned}
\eqno(\NUM)$$  %(3.13)$$
Observe that (\theIIIpxiii.\theiiipxiii(i))\   %(3.13(i)) 
holds automatically with the trivial value ${\bf 0}$.
Furthermore conditions 
(\theIIIpxiii.\theiiipxiii(i)) and (\theIIIpxiii.\theiiipxiii(ii))\   %(3.13(i)) and (3.13(ii)) 
are equivalent
(by changing the roles of ${\bf a}$ and ${\bf b}$).
Finally we observe that, in (\theIIIpxiii.iiipxiii(i)),  %(3.13(i)), 
for fixed 
${\bf a,p}$ and ${\bf y}\in{\rm Co}\{ {\bf a,p}\}$
we can choose the points ${\bf b}$ and $\widetilde{\bf b}$ 
on different half plain components 
of $\RR^2 \!\setminus\!{\rm Line}\{ {\bf a,p}\}$ arbitrarily.
This implies that all the vectors 
$\nabla f^{i,{\bf b}}_{{\bf a,p}} ({\bf y})$,
$\nabla f^{i, \widetilde{\bf b}}_{{\bf a,p}} ({\bf y})$
with ${\bf b},\widetilde{\bf b}\in \RR^2 \!\setminus\!{\rm Line}\{ {\bf a,p}\}$
must be the same.
Due to the construction (\theIpiv.\theipiv), %(1.4), 
the fact that all the pairs 
$\varphi_{\bf p,T},\varphi_{\bf p,\widetilde{T}}$ 
resp.
$\psi^{(j)}_{\bf p,T},\psi^{(j)}_{\bf p,\widetilde{T}}$
of basic functions 
touch ${\cal C}^1$-smoothly in case of adjacent triangles
${\bf T,\widetilde{T}}$, ensures that 
the splines $f_{{\cal T},F}$ are all ${\cal C}^1$-smooth as well.    
 
\ADT
\newcounter{IIIpvi} \setcounter{IIIpvi}{\value{SEC}}
\newcounter{iiipvi} \setcounter{iiipvi}{\value{THM}}
\bigskip
\noindent
{\bf Notation \NUM.}  %3.6.} 
Given any ordered triple
$({\bf u},{\bf v},{\bf w})$ % \in \big[\RR^2\big]^3$ of
of non-collinear points, we shall write
${\bf g}_{\bf u,v}^{\bf w} :=\nabla\lambda_{\bf u,v}^{\bf w}$  
for the  
constant gradient vectors of the baricentric coordinate functions.
Notice that, by $(\theIIpi.\theiipi)$,  %$(2.1)$,
\ADT
\newcounter{IIIpx} \setcounter{IIIpx}{\value{SEC}}
\newcounter{iiipx} \setcounter{iiipx}{\value{THM}}
$${\bf g}_{\bf u,v}^{\bf w} := 
\frac{ {\bf (u-v)R} }{ \big\langle {\bf (u-v)R} \big\vert {\bf w-u}\big\rangle}
= \frac{ \sigma_{\bf u,v}^{\bf w} {\bf (u-v)R} }{ 
{\rm area}({\rm Co}\{ {\bf u,v,w}\}) } . 
\eqno(\NUM)$$  %(3.10)$$
where \ $\sigma_{\bf u,v}^{\bf w} \!=\! \pm 1$ \ according as 
\ $({\bf u,v,w})$ \ are oriented anticlockwise or clockwise.
In particular, if \ ${\bf T} = {\rm Co}\{ {\bf a,b,p}\}$ \  is a
non-degenerate triangle, we have
%%\ADT
%%\newcounter{IIIpix} \setcounter{IIIpix}{\value{SEC}}
%%\newcounter{iiipix} \setcounter{iiipix}{\value{THM}}
$$\begin{aligned}
&{\bf g}_{\bf a,b}^{\bf p} + {\bf g}_{\bf b,p}^{\bf a} +
{\bf g}_{\bf a,p}^{\bf b} = 
\nabla \big[ \lambda_{\bf a,b}^{\bf p} + \lambda_{\bf b,p}^{\bf a} +
\lambda_{\bf a,p}^{\bf b} \big] = \nabla 1 = {\bf 0} ,\\
&{\bf g}_{\bf a,b}^{\bf p} \perp {\bf b-a},\quad 
{\bf g}_{\bf a,p}^{\bf b} \perp {\bf a-p}, \quad 
{\bf g}_{\bf b,p}^{\bf a}\perp {\bf b-p} .
\end{aligned}$$
%%\eqno(\NUM)$$  %(3.9)$$

\ADT
\newcounter{IIIpxix} \setcounter{IIIpxix}{\value{SEC}}
\newcounter{iiipxix} \setcounter{iiipxix}{\value{THM}}
\bigskip
\noindent
{\bf Lemma \NUM.}  %3.19.} 
\it If ${\bf T}={\rm Co}\{ {\bf a,b,p}\}$ 
is a non-degenerate triangle, at the points
${\bf y}_t := (1-t){\bf a} + t{\bf p}$ of the edge 
${\rm Co}\{ {\bf a,p}\}$ we have
\ADT
\newcounter{IIIpxixw} \setcounter{IIIpxixw}{\value{SEC}}
\newcounter{iiipxixw} \setcounter{iiipxixw}{\value{THM}}
$$\begin{aligned}
\nabla f^{i,{\bf p}}_{\bf a,b} ({\bf y}_t) =
&x^{[i]}\big((1\!-\! t)({\bf a \!-\! p}) \big)   
[\Phi^{[i]}]^\prime \! (t) \, {\bf g}^{\bf p}_{\bf a,b} +\\ 
&+ \Phi^{[i]}(t) {\bf e}^{[i]} +
t^2 (1\!- \!t) P^{i,{\bf p}}_{\bf a,b}(0,1\!- \! t) \, 
{\bf g}^{\bf b}_{\bf a,p} .
\end{aligned}
\eqno(\NUM)$$   %(3.19^\prime)$$

\noindent
{\bf Proof.} \rm With the abbreviations 
$$\ell_0:=\lambda_{\bf a,b}^{\bf p}, \ 
\ell_1:=\lambda_{\bf a,p}^{\bf b}, \ \ell_2:=\lambda_{\bf b,p}^{\bf a}, \quad
P^{[i]}:= P^{i,{\bf p}}_{{\bf a},{\bf b}} , \quad
G^{[i]} := \ell_0^2 \ell_2 P^{[i]} (\ell_1,\ell_2)$$ 
we can write
%%\ADT
%%\newcounter{IIIpxx} \setcounter{IIIpxx}{\value{SEC}}
%%\newcounter{iiipxx} \setcounter{iiipxx}{\value{THM}} 
$$\begin{aligned}
\nabla f_{{\bf a},{\bf b}}^{i,{\bf p}} &=
\nabla\Big[ x^{[i]}_{\bf p} \Phi^{[i]}(\ell_0) + 
\ell_1  G^{[i]} \Big] =\\
&= x^{[i]}_{{\bf p}} \,\nabla\! \big[ {\Phi^{[i]}}(\ell_0) \big] +
{\Phi^{[i]}}(\ell_0) \nabla x^{[i]}_{{\bf p}} + \ell_1 \nabla G^{[i]} +
G^{[i]} \nabla \ell_1 = \\
&= x^{[i]}_{{\bf p}} {\Phi^{[i]}}^\prime(\ell_0) \nabla\ell_0 +
\Phi^{[i]}(\ell_0) {\bf e}^{[i]} + 
\ell_1 \nabla G^{[i]} + G^{[i]}  \nabla\ell_1 .
\end{aligned}$$
%%\eqno(\NUM)$$  %(3.20)$$
We complete the proof with the observations that
%%\ADT
%%\newcounter{IIIpxxi} \setcounter{IIIpxxi}{\value{SEC}}
%%\newcounter{iiipxxi} \setcounter{iiipxxi}{\value{THM}}
$$\begin{aligned}
&\ell_0({\bf y}_t) \!=\! t,\quad \ell_1({\bf y}_t) \!=\! 0,\quad 
\ell_2({\bf y}_t) \!=\! 1\!-\! t, \\
&x^{[i]}_{\bf p} ({\bf y}_t) \!=\! 
x^{[i]} \big( (1\!- \! t)({\bf a\!-p}) \big) , \quad 
\nabla x^{[i]}_{\bf p} \!\equiv\! {\bf e}^{[i]} .
\end{aligned}$$
%%\eqno(NUM)$$  %(3.21)$$

\ADT
\newcounter{IIIpxxii} \setcounter{IIIpxxii}{\value{SEC}}
\newcounter{iiipxxii} \setcounter{iiipxxii}{\value{THM}}
\noindent
{\bf Remark \NUM.}  %3.22.}  
To prove Theorem \theIIpvi.\theiipvi,  %2.6,
we need a precise description for the coefficients of the polynomials
$P^{i,{\bf p}}_{\bf a,b}$ in terms of the variables ${\bf a,b,p}$
such that $(\theIIIpxvw.\theiiipxvw)$  %$(3.15^\prime)$ 
should hold. 

\medskip
\rm
According to Lemma \theIIIpxv.\theiiipxv,  %3.15, 
the procedure $\SS : ({\cal T},F)\!\mapsto\! f_{{\cal T},F}$
produces ${\cal C}^1$-splines for every admissible data
if and only if, for any $t\!\in\![0,1]$ and 
for any fixed pair ${\bf a,p}$ of distinct points,
the gradient expressions $(\theIIIpxixw.\theiiipxixw)$\  %$(3.19^\prime)$
are independent of the variable ${\bf b}$ ranging in 
$\RR^2 \!\setminus\!{\rm Line}\{ {\bf a,p}\}$.
This latter 
condition can be formulated
in terms of the 
${\bf b}$-independent affine coordinates $(\theIIpiw.\theiipiw)$  %$(2.1^\prime)$ 
as follows.
By $(\theIIIpx.\theiiipx)$\  %$(3.10)$ 
we have
$$\begin{aligned}
{\bf g}^{\bf b}_{\bf a,p} \!&=\! 
\frac{{\bf (a\!-\! p)R}}{\langle {\bf (a\!-\! p)R} 
\vert {\bf b\!-\! a}\rangle} \!=\!
\Vert {\bf p \!-\! a}\Vert^{-2} 
(1/\overline{\xi}^{\ {\bf b}}_{\bf p,a}) {\bf (p \!-\!a)R} ,\\
{\bf g}^{\bf p}_{\bf a,b} \!&=\!
\frac{{\bf (a\!-\! b)R}}{\langle {\bf (a\!-\! b)R} 
\vert {\bf p\!-\! a}\rangle} \!=\!
\frac{\xi^{\bf (a\!-\! b)R}_{\bf p,a} {\bf (p\!-\! a)} \!+\!
\overline{\xi}^{\bf (a\!-\! b)R}_{\bf p,a} {\bf (p\!-\! a)R}}{\Vert 
{\bf p\!-\! a}\Vert^2 \overline{\xi}^{\bf b}_{\bf p,a}} =\\
&= \Vert {\bf p \!-\! a}\Vert^{-2} \Big[ {\bf (p \!-\! a)} \!+\! 
(\xi^{\bf b}_{\bf p,a}/\overline{\xi}^{\ {\bf b}}_{\bf p,a}) 
{\bf (p \!-\! a)R} \Big] . 
\end{aligned}
$$
Thus we can rewrite $(\theIIIpxixw.\theiiipxixw)$\  %$(3.19^\prime)$ 
in the form
\ADT
\newcounter{IIIpxixww} \setcounter{IIIpxixww}{\value{SEC}}
\newcounter{iiipxixww} \setcounter{iiipxixww}{\value{THM}}
$$\begin{aligned}
&\nabla f^{i,{\bf p}}_{\bf a,b}({\bf y}_t) =
\Big[{\bf b}\hbox{-independent terms}\Big] + \\
&+ \frac{ x^{[i]}\big( (1 \!-\! t) ({\bf a\!-\! p}) \big) 
[\Phi^{[i]}]^\prime \!(t) 
\xi^{\bf b}_{\bf p,a} +
t^2 (1\!-t) P^{i,{\bf p}}_{\bf a,b}(0,1\!-t) 
}{ \Vert {\bf p \!-\! a} \Vert^2 
\overline{\xi}^{\ {\bf b}}_{\bf p,a} } 
{\bf (p \!-\! a)R} .
\end{aligned}
\eqno(\NUM)$$   %(3.19^{\prime\prime})$$
Hence we conclude immediately the following.

\ADT
\newcounter{IIIpxxiv} \setcounter{IIIpxxiv}{\value{SEC}}
\newcounter{iiipxxiv} \setcounter{iiipxxiv}{\value{THM}}
\bigskip
\noindent{\bf Lemma \NUM.}  %3.24.} 
\it 
We have $(\theIIIpxvw.\theiiipxvw)$  %$(3.15^\prime)$
if and only if
for every pair ${\bf p,a}$ of distinct points 
there exist polynomials
$K^{i,{\bf p}}_{\bf a}$ $(i\!=\! 0,1,2)$ 
of one variable
such that
\DEL{
\begin{aligned}
K^{i,{\bf p}}_{\bf a}(t) \equiv
\big( \overline{\xi}^{\ {\bf b}}_{\bf p,a} \big)^{-1} \Big[
&x^{[i]}\big( (1 \!-t) ({\bf a\!-\! p})\big) [\Phi^{[i]}]^\prime \! (t) 
\xi^{\bf b}_{\bf p,a} +\\ 
&+ t^2 (1\!-t) P^{i,{\bf p}}_{\bf a,b}(0,1\!-t) \Big] 
%\equiv K^{i,{\bf p}}_{\bf a}(t)
\end{aligned}
}%ENDDEL
\ADT
\newcounter{IIIpxxiii} \setcounter{IIIpxxiii}{\value{SEC}}
\newcounter{iiipxxiii} \setcounter{iiipxxiii}{\value{THM}}
$$K^{i,{\bf p}}_{\bf a}(t) =
x^{[i]}\big( (1 \!-t) ({\bf a\!-\! p})\big) [\Phi^{[i]}]^\prime \! (t) 
\frac{\xi^{\bf b}_{\bf p,a}}{\overline{\xi}^{\ {\bf b}}_{\bf p,a}}  + 
\frac{t^2 (1\!-t)}{\overline{\xi}^{\ {\bf b}}_{\bf p,a}} 
P^{i,{\bf p}}_{\bf a,b}(0,1\!-t)
\eqno(\NUM)$$  %(3.23)$$
independently of the choice of ${\bf b}$ 
outside ${\rm Line}\{ {\bf a,p}\}$.

\bigskip
\rm 
We can regard $(\theIIIpxxiii.\theiiipxxiii)$\    %$(3.23)$ 
as a partial algebraic condition 
on the polynomials $P^{i,{\bf p}}_{\bf a,b}$ of two variables as
\ADT
\newcounter{IIIpxxiiiw} \setcounter{IIIpxxiiiw}{\value{SEC}}
\newcounter{iiipxxiiiw} \setcounter{iiipxxiiiw}{\value{THM}} 
$$\begin{aligned}
P^{i,{\bf p}}_{\bf a,b}(0,1\!-t) =\  &\overline{\xi}^{\ {\bf b}}_{\bf p,a} 
\frac{ K^{i,{\bf p}}_{\bf a}(t) }{ t^2 (1\!-\! t)} -\\ 
&- \xi^{\ {\bf b}}_{\bf p,a} 
\frac{ x^{[i]}\big( (1\!-\! t)({\bf a\!-\! p})\big) [\Phi^{[i]}]^\prime \!(t) 
}{ t^2(1\!-\! t)} \qquad (0 \!<\! t \!<\! 1).
\end{aligned}
\eqno(\NUM)$$  %(3.23^\prime)$$
Since, for fixed ${\bf a,p}$, the coordinates  
$\big( \xi^{\ {\bf b}}_{\bf p,a} , \overline{\xi}^{\ {\bf b}}_{\bf p,a} \big)$
may assume arbitrary values $(r,s)$ with $s\ne 0$,  
from $(\theIIIpxxiiiw.\theiiipxxiiiw)$\  %$(3.23^\prime)$ 
we obtain the polynomial divisability relations \ 
$t^2(1-t) \big\vert K^{i,{\bf p}}_{\bf a}(t)$ 
and
$t^2(1-t) \big\vert 
x^{[i]}\big((1\!-\! t)({\bf a \!-\! p}) \big)  [\Phi^{[i]}]^\prime \! (t)$.
Since $x^{[0]}\big( (1-t){\bf (a-p)}\big) \equiv 1$ 
and
$x^{[0]}\big( (1-t){\bf (a-p)}\big) \equiv (1-t)x^{[j]}({\bf a-p})$
for $j=1,2$,
with the aid of $(3.22^\prime)$
we can state $(\theIIIpxxiiiw.\theiiipxxiiiw)$\  %$(3.23^\prime)$ 
in the form
\ADT
\newcounter{IIIpxxvi} \setcounter{IIIpxxvi}{\value{SEC}}
\newcounter{iiipxxvi} \setcounter{iiipxxvi}{\value{THM}}
$$P^{i,{\bf p}}_{{\bf a},{\bf b}} (0,1 \!-\! t) =
\frac{ \big\langle {\bf b \!-\! a} \big\vert {\bf p \!-\! a} \big\rangle
}{ \Vert {\bf p \!-\! \bf a}\Vert^2 }\ x^{[i]}_{\bf p}({\bf a})\  
\chi^{[i]}(t) +
\frac{ \big\langle {\bf b \!-\! a} \big\vert ({\bf p \!-\! a}){\bf R} \big\rangle
}{ \Vert {\bf p \!-\! a}\Vert^2 }\ 
\kappa^{i,{\bf p}}_{\bf a}(t) 
\eqno(\NUM)$$  %(3.26)$$
with suitable polynomials $\chi^{[i]}$ and 
$\kappa^{i,{\bf p}}_{\bf a}$ $(i=0,1,2;\ {\bf a \ne p}\in\RR^2)$
of one variable.
Actually
$$\begin{aligned}
&\kappa^{i,{\bf p}}_{\bf a}(t) = \frac{K^{i,{\bf p}}_{\bf a}(t) }{ t^2(1-t)^2}, 
\quad
\chi^{[0]}(t) = \frac{[\Phi^{[0]}]^\prime(t) }{ t^2(1-t)} = 
\frac{\Phi^\prime(t) }{ t^2(1-t)} , \\ 
&x^{[i]}_{\bf p}({\bf a})\chi^{[j]}(t) \!=\! 
\frac{x^{[j]}\big( (1 \!-\! t){\bf (a \!-\! p)}\big)[\Phi^{[j]}]^\prime(t) }{ t^2(1-t)} 
%= x^{[j]}_{\bf p}({\bf a})  
%\frac{\Theta^\prime(t) }{ t^2} \!=\\ 
%&\hskip20mm 
= x^{[j]}_{\bf p}({\bf a}) \frac{[\Psi(t)/(t \!-\!1)]^\prime }{ t^2}
\end{aligned}$$
for $j\!\!=\!\! 1,2$ on the basis of $(\theIIIpxxiiiw.\theiiipxxiiiw)$\   %$(3.23^\prime)$.
In terms of the Kronecker-$\delta$, 
we can write even
$$\chi^{[i]}(t) = t^{-2} (1-t)^{-\delta_{i,0}} 
[\Phi^{[i]}]^\prime(t) \qquad (i=0,1,2).$$ 
Clearly, the polynomials $K^{i,{\bf p}}_{\bf c}$ cannot be chosen
arbitrarily. There is a unique obstacle:
we obtained Lemma \theIIIpxix.\theiiipxix\   %3.19 
and hence $(\theIIIpxxiii.\theiiipxxiii)$\   %$(3.23)$ 
by an inspection
of $\nabla f^{i,{\bf p}}_{\bf a,b}$ on one of the edges
of a triangle ${\bf T}={\rm Co}\{ {\bf a,b,p}\}$ at the 
distinguished point ${\bf p}$ 
(namely ${\rm Co}\{ {\bf a,p}\}$ with the parametrization
${\bf y}_t := (1-t){\bf a} + t{\bf b}$)
while also the analogous conclusion should also be taken 
simultaneously in to account with the second edge 
(namely ${\rm Co}\{ {\bf b,p}\}$ issued from ${\bf p}$.
Applying a change ${\bf a}\leftrightarrow{\bf b}$ 
and taking into account the symmetry 
$(\theIIIpxviiww.\theiiipxviiww)$,  %$(3.17^{\prime\prime})$, 
we see that also
\ADT
\newcounter{IIIpxxviw} \setcounter{IIIpxxviw}{\value{SEC}}
\newcounter{iiipxxviw} \setcounter{iiipxxviw}{\value{THM}} 
$$P^{i,{\bf p}}_{{\bf a},{\bf b}} (1 \!-\! t,0) \!=\!
\frac{ \big\langle {\bf a \!-\! b} \big\vert {\bf p\!-\! b} \big\rangle
}{ \Vert {\bf p \!-\! b}\Vert^2 } x^{[i]}_{\bf p}({\bf b})\ 
\chi^{[i]}(t) \!+\!
\frac{ \big\langle {\bf a\!-\! b} \big\vert ({\bf p\!-\! b}){\bf R} \big\rangle
}{ \Vert {\bf p\!-\! b}\Vert^2 }\ 
\kappa^{i,{\bf p}}_{\bf b}(t) .
\eqno(\NUM)$$  %(3.26^\prime)$$
We obtain the complete description for the families of 
polynomials $K^{i,{\bf p}}_{\bf a,b}$ being admissible 
by Lemma \theIIIpxxiv.\theiiipxxiv\  %3.24
by the next obervation.

\ADT
\newcounter{IIIpxxvii} \setcounter{IIIpxxvii}{\value{SEC}}
\newcounter{iiipxxvii} \setcounter{iiipxxvii}{\value{THM}}
\bigskip
\noindent{\bf Lemma \NUM.}  %3.27.} 
\it For any pair 
${\bf p \ne c}\in\RR^2$, in $(\theIIIpxxvi.\theiiipxxvi\!-\!\theiiipxxviw)$\ %$(3.26-26^\prime)$ 
we have \ 
$\chi^{[i]}(1) = 
\kappa^{i,{\bf p}}_{\bf c}(1) =0$.
 
\bigskip
\noindent
{\bf Proof}. \rm Fix $i,{\bf p} \in\RR^2$ and $\rho>0$ arbitrarily. 
Consider 
$(\theIIIpxxvi.\theiiipxxvi\!-\!\theiiipxxviw)$\   %$(3.26-26^\prime)$ 
for pairs
${\bf a,b}$ with $\Vert {\bf a-p}\Vert = \Vert {\bf b-p}\Vert =\rho$
written in the form   
$${\bf a := c}_\sigma, \ \ {\bf b := c}_\tau \quad {\rm where} \quad 
{\bf c}_\tau := {\bf p} + \rho {\bf u}_\tau, \ 
{\bf u}_\tau :=  \cos\tau {\bf e}^{[1]} + \sin\tau {\bf e}^{[2]}  .$$
Due to  
$(\theIIIpxviiww.\theiiipxviiww)$,  %$(3.17^{\prime\prime})$,
with the abbreviations  $\alpha:= \chi^{[i]}(1)$ and 
$\beta(\tau) := \kappa^{i,{\bf p}}_{{\bf c}_\tau}(1)$ 
we get
$$\begin{aligned}
&0 = P^{i,{\bf p}}_{{\bf c}_\sigma,{\bf c}_\tau}(0,0) -
P^{i,{\bf p}}_{{\bf c}_\tau,{\bf c}_\sigma}(0,0) =\\
&= \Big[ 
\big( \langle {\bf u}_\tau \vert {\bf u}_\sigma\rangle \!-\! 1\big)
x^{[i]} ({\bf c}_\sigma) \alpha 
\!+\! \langle {\bf u}_\tau \vert {\bf u}_\sigma {\bf R}\rangle \beta(\sigma)\Big] -\\
&\hskip6mm - \Big[ \big( \langle {\bf u}_\sigma \vert {\bf u}_\tau\rangle \!-\! 1\big) 
x^{[i]} ({\bf c}_\tau)\alpha 
\!+\! \langle {\bf u}_\sigma \vert {\bf u}_\tau {\bf R}\rangle
\beta(\tau) \Big]  \!=\\
&= ( \langle {\bf u}_\sigma \vert {\bf u}_\tau\rangle \!-\! 1 )
\big[ x^{[i]} ({\bf c}_\sigma) - x^{[i]} ({\bf c}_\tau) \big] \alpha +
\langle {\bf u}_\tau \vert {\bf u}_\sigma {\bf R}\rangle \beta(\sigma) - \langle {\bf u}_\sigma \vert {\bf u}_\tau {\bf R}\rangle \beta(\tau) \big) =\\
&= \big[ \cos(\sigma-\tau) - 1 \big] 
\big[ x^{[i]} ({\bf c}_\sigma) - x^{[i]} ({\bf c}_\tau) \big] \alpha + 
\sin(\tau-\sigma) \big[ \beta(\sigma)+\beta(\tau) \big] .
\end{aligned} 
$$
Since $x^{[0]} \equiv 1$, in any case we have 
$x^{[i]}_{\bf p}({\bf c}_\sigma) - x^{[i]}_{\bf p}({\bf c}_\tau) =
\rho [ x^{[i]} ({\bf u}_\sigma) -  x^{[i]} ({\bf u}_\tau) ]$.
It follows
\ADT
\newcounter{IIIpxxix} \setcounter{IIIpxxix}{\value{SEC}}
\newcounter{iiipxxix} \setcounter{iiipxxix}{\value{THM}}
$$\begin{aligned}
\beta(\sigma) + \beta(\tau) &= \alpha \rho 
\frac{\cos(\tau-\sigma)-1 }{ \sin(\tau-\sigma) } 
\big[ x^{[i]} ({\bf u}_\sigma) -  x^{[i]} ({\bf u}_\tau) \big] ,\\ 
\big\vert \beta(\sigma) + \beta(\tau)  \big\vert &\le 
\rho \vert\alpha \vert  
\frac{1 -\cos(\tau-\sigma) }{ \sin( \vert\tau-\sigma\vert ) } 
\big\Vert {\bf u}_\sigma - {\bf u}_\tau \big\Vert \le\\
&\le 2 \rho \vert\alpha \vert [ 1 -\cos(\tau-\sigma) ] .
\end{aligned}
\eqno(\NUM)$$  %(3.29)$$
Suppose indirectly $\beta(\tau) \ne 0$ for some $\tau\in\RR$. 
Let 
$\eps := \vert \beta(\tau)\vert$ and choose $\delta> 0$
such that
$2 \rho \vert\alpha \vert ( 1 -\cos \theta ) <\eps /4$
whenever $\vert\theta\vert \le \eps$.
Then we have  
$\vert \beta(\tau)+ \beta(\tau \pm \delta/2)\vert < \eps/4$
that is 
$\beta(\tau \pm\delta/2) \in \big[ -\eps/4,\eps/4] - \beta(\tau)$.
Therefore
$\beta(\tau+\delta/2)+ \beta(\tau - \delta/2) \in 
\big[ -\eps/2,\eps/2 \big] - 2\beta(\tau)
\subset \big[ -\eps/2,\eps/2 \big] + \{ -2\eps, 2\eps\} =
\big[ -5\eps/2,-3\eps/2 \big] \cup \big[ 3\eps/2,5\eps/2 \big]$
i.e. 
$\vert \beta(\tau+\delta/2)+ \beta(\tau - \delta/2)\vert \in 
\big[ 3\eps/2,5\eps/2 \big]$
However, we also have
$\vert \beta(\tau+\delta/2)+ \beta(\tau - \delta/2)\vert < \eps/4$
which leads to the contradiction
$\vert \beta(\tau+\delta/2)+ \beta(\tau - \delta/2) \vert \in 
\big[ 3\eps/2, 5\eps/2 \big] \cap \big[ 0,\eps/4 \big] =\emptyset$.
By the arbitrariness of 
the radius $\rho$, the angle $\tau$ and the origin ${\bf p}$,
we conclude that $\kappa^{i,{\bf p}}_{\bf c}(1)=0$ in any case. 

For $i=1,2$ we get $\alpha=0$ i.e. $\chi^{[i]}(1)=0$ immediately
by plugging $\beta(\tau)=\beta(\sigma)= 0$ with $\sigma:=\tau+\pi/4$ in 
the first equation of $(\theIIIpxxix.\theiiipxxix)$.  $(3.29)$.
(Remark:
$x^{[0]}({\bf u}_\sigma) - x^{[0]}({\bf u}_\tau) =1-1=0$,
thus the argument does not work for $i=0$).
In the case $i=0$ we conclude $\alpha=0$ as follows.
Consider the difference of equations 
$(\theIIIpxxvi.\theiiipxxvi\!-\!\theiiipxxviw)$\   %$(3.26-26^\prime)$
for $t=1$ with 
${\bf a}:={\bf p} + {\bf e}^{[1]}$ and
${\bf b}:={\bf p} + {\bf e}^{[1]} + {\bf e}^{[2]}$.
Since $\kappa^{i,{\bf p}}_{\bf c} =0$ $({\bf c}={\bf a,b}$
is estabished already, we get simply
$0 = -(1/2)\chi^{[0]}(1)$ which completes the proof.

\ADT
\newcounter{IIIpxxviii} \setcounter{IIIpxxviii}{\value{SEC}}
\newcounter{iiipxxviii} \setcounter{iiipxxviii}{\value{THM}}
\bigskip
\noindent
{\bf Corollary \NUM.}  %3.28.} 
\it The relations $(\theIIIpxvw.\theiiipxvw)$  %$(3.15^\prime)$ 
hold 
if and only if we have $(\theIIIpxxiii.\theiiipxxiii)$\  %$(3.23)$ 
with the symmetry $(\theIIIpxviiww.\theiiipxviiww)$  %$(3.17^{\prime\prime})$ 
where the polynomials \  
$K^{i,{\bf p}}_{\bf c}(t)$ \  respectively 
$x^{[i]}\big((1\!-\! t){\bf (a-p)}\big) [\Phi^{[i]}]^\prime(t)$ \ 
are all divisable by \ $t^2(1-t)^2$. 
 
\bigskip
\noindent
{\bf Proof.} \rm
The relation $\kappa^{i,{\bf p}}_{\bf c}(1) \!=\! 0$ implies that 
there is a polynomial  $\widetilde{\kappa}^{i,{\bf p}}_{\bf c}$
such that
$\kappa^{i,{\bf p}}_{\bf c}(t) = 
(1 \!-\! t) \widetilde{\kappa}^{i,{\bf p}}_{\bf c}(t)$
and
$ K^{i,{\bf p}}_{\bf c}(t) = 
t^2(1 \!-\! t) \kappa^{i,{\bf p}}_{\bf c} = 
t^2(1 \!-\! t)^2 \widetilde{\kappa}^{i,{\bf p}}_{\bf c}(t)$
with some polynomial.
Similarly, from $\chi^{[i]}(1)=0$ we conclude that
$\chi^{[i]}(t) = (1 \!-\! t) \widetilde{\chi}^{[i]}(t)$
and 
$(1\!-\! t) x^{[i]}({\bf a}) [\Phi^{[i]}]^\prime(t) = 
((1\!-\! t) t^2 \chi^{[i]} (t) = 
t^2(1 \!-\! t)^2 \widetilde{\chi}^{[i]}(t)$
with some polynomial $\widetilde{\chi}^{[i]}$.

\ADT
\newcounter{IIIpxxviiiw} \setcounter{IIIpxxviiiw}{\value{SEC}}
\newcounter{iiipxxviiiw} \setcounter{iiipxxviiiw}{\value{THM}}
\bigskip
\noindent
{\bf Corollary \NUM.}  %3.28$^\prime$.} 
\it We can write
$K^{i,{\bf p}}_{\bf c}(t) = t^2(1-t)^2 k^{i,{\bf p}}_{\bf c}(t)$
$({\bf p\ne c}\in\RR^2)$ and the admissible shape functions
$\Phi,\Psi$ have the form
\ADT
\newcounter{IIIpxxxi} \setcounter{IIIpxxxi}{\value{SEC}}
\newcounter{iiipxxxi} \setcounter{iiipxxxi}{\value{THM}}
$$\begin{aligned}
&{\rm (i)} \ \ \ \Phi(t) = t^3(10-15t+6t^2)+t^3(1-t)^3 \Phi_1(t),\\ 
&{\rm (ii)} \ \ \Psi(t) = t^3(t-1)(4-3t) + t^3(1-t)^3\Psi_1(t) 
\end{aligned} 
\eqno(\NUM)$$  %(3.31)$$
with suitable polynomials $k^{i,{\bf p}}_{\bf c},\Phi_2,\Psi_2$.

\bigskip
\noindent
{\bf Proof.} \rm The stated form of $K^{i,{\bf p}}_{\bf c}$ is
clear from $\theIIIpxxviii.\theiiipxxviii)$.   %$3.28$. 
By definition
$\Phi^{[0]}(t) = \Phi(t)$ and
$x^{[0]}\big( (1-t){\bf (a-p)}\big) \equiv 1$.
Furthermore
$\Phi^{[j]}(t) = \Psi(t)/(t-1)$ 
and $x^{[j]}\big( (1 \!-\! t){\bf (a \!-\! p)}\big) \!=\!
(1 \!-\! t)x^{[j]}_{\bf p}({\bf a})$
for $j \!=\! 1,2$.
Thus, taking $(\theIIIpi.\theiiipi)$  %$(3.1)$ 
into acount,
the relation that \ $t^2(1 \!-\! t)^2$ \ is a divisor of \ 
$x^{[0]}\big( (1 \!-\! t){\bf (a \!-\! p)}\big) [\Phi^{[0]}]^\prime(t) \!=\! 
\Phi^\prime(t) \!=\!
6t(1 \!-\! t) + 2t(1 \!-\! t)(1\!-\! 2t)\Phi_0(t) + t^2(1 \!-\! t)^2\Phi^\prime_0(t)$
means simply that 
$t(1\!-\! t) \big\vert 6+2(1\!-\! 2t)\Phi_0(t)$ 
i.e. 
$6 + 2(1-2t)\Phi_0(t)\vert_{t=0,1} =0$ implying
$\Phi_0(0)=-3$, $\Phi_0(1)=3$.
Therefore
$\Phi_0(t) = -3+6t + t(1-t)\Phi_1(t)$
with a polynomial $\Phi_1$ and
the generic form of $\Phi$ is
$\big(\theIIIpxxxi.\theiiipi{\rm (i)}\big)$.   %$\big(3.31{\rm (i)}\big)$
Also according to $(\theIIIpi.\theiiipi)$  %$(3.1)$, 
in the cases $j=1,2$ we can write
$\Psi(t)=-t^2(1-t)+t^2(1-t)^2\Psi_0(t)$ with some polynomial $\Psi_0$.
Thus the relation that
$t^2(1-t)^2$ is a divisor of 
$x^{[j]}\big( (1-t){\bf (a-p)}\big) [\Phi^{[0]}]^\prime(t) \equiv
(1-t) x^{[j]}_{\bf p}({\bf a}) \big[\Psi(t) /(1-t) \big]^\prime$
means that
$t^2(1-t) \big\vert \big[\Psi(t) /(1-t) \big]^\prime \equiv 
-2t + t(2-3t)\Psi_0(t) + t^2(1-t)^2\Psi_0^\prime(t)$
is equivalent to saying
$t(1-t)\big\vert -2 + (2-3t)\Psi_0(t)$
i.e.
$-2 + (2-3t)\Psi_0(t) \vert_{t=0,1} =0$
implying $\Psi_0(0)=1$ and $\Psi_0(1)=-2$.
Therefore
$\Psi_0(t) = 1 - 3t + t^2(1-t)\Psi_1(t)$
with some polynomial $\Psi_1$ and
the generic form of $\Psi$ is
$\big(\theIIIpxxxi.\theiiipxxxi{\rm (i)}\big)$.   %$\big(3.31{\rm (ii)}\big)$.

\ADT
\newcounter{IIIpxxxix} \setcounter{IIIpxxxix}{\value{SEC}}
\newcounter{iiipxxxix} \setcounter{iiipxxxix}{\value{THM}}
\bigskip
\noindent
{\bf \NUM. Finish of the proof of Theorem \theIIpvi.\theiipvi}  %2.6}

\medskip 

In view of $(\theIIIpxxvi.\theiiipxxviw\!-\!\theiiipxxviw)$\   %$(3.26^\prime-26^{\prime\prime})$ 
and Remark \theIIIpii.\theiiipii(ii)\   %3.2(ii) 
we can write
%%\ADT
%%\newcounter{IIIpxxx} \setcounter{IIIpxxx}{\value{SEC}}
%%\newcounter{iiipxxx} \setcounter{iiipxxx}{\value{THM}}
$$\begin{aligned}
&P^{i,{\bf p}}_{\bf a,b}(s,t) =  
P^{i,{\bf p}}_{\bf a,b}(0,0) + 
s \big[ \big( P^{i,{\bf p}}_{\bf a,b}(s,0) - 
P^{i,{\bf p}}_{\bf a,b}(0,0) \big)/s \big] +\\ 
&\hskip 20mm + t \big[ \big( P^{i,{\bf p}}_{\bf a,b}(0,t) - 
P^{i,{\bf p}}_{\bf a,b}(0,0) \big)/t \big] + 
st {\rm Pol}(s,t) =\\
&= s \big[ \big( P^{i,{\bf p}}_{\bf b,a}(0,s) /s \big] +
t \big[ \big( P^{i,{\bf p}}_{\bf a,b}(0,t)/t \big] + st {\rm Pol}(s,t) =\cr
&= s \left[ \overline{\xi}^{\ {\bf a}}_{\bf p,b} 
\frac{ K^{i,{\bf p}}_{\bf b}(1\!-\! s) }{ s^2 (1\!-s)^2} - 
\xi^{\ {\bf a}}_{\bf p,b} 
\frac{ x^{[i]}({\bf b}) \big[\Phi^{[i]}\big]^\prime \!(1\!-\! s) 
}{ s (1\!-\! s)^2} \right] +\\
&\qquad +  t \left[ \overline{\xi}^{\ {\bf b}}_{\bf p,a} 
\frac{ K^{i,{\bf p}}_{\bf a}(1\!-\! t) }{ t^2 (1\!-t)^2} - 
\xi^{\ {\bf b}}_{\bf p,a} 
\frac{ x^{[i]}({\bf a}) \big[\Phi^{[i]}\big]^\prime \!(1\!-\! t) 
}{ t (1\!-\! t)^2} \right] + st R^{\bf p}_{\bf a,b}(s,t) 
\end{aligned}$$
%%\eqno(\NUM)$$  %(3.30)$$
with suitable  
polynomials 
$K^{i,{\bf p}}_{\bf c}, \Phi^{[i]}, R^{\bf p}_{\bf a,b}$ 
of one- resp. two variables
such that
$t^2(1 \!-\! t)^2 \big\vert K^{i,{\bf p}}_{\bf c}(t)$ 
and
$t^2(1 \!-\! t) \big\vert x^{[i]}_{\bf p}({\bf a}) [\Phi^{[i]}]^\prime(t)$.
It is straightforward to check that
the functions $f^{i,{\bf p}}_{\bf a,b}$ are polynomials in these cases 
and 
$P^{i,{\bf p}}_{\bf a,b} (s,t) = P^{i,{\bf p}}_{\bf b,a} (t,s)$
if and only if 
$R^{\bf p}_{\bf a,b} (s,t) = R^{\bf p}_{\bf b,a} (t,s)$.
It remains to show that the expressions
$$\begin{aligned}
&\nabla f^{i,{\bf p}}_{\bf a,b}({\bf y}_t) \quad  
{\rm with} \quad {\bf y}_t:= (1-t){\bf a}+t{\bf p} , \\ 
%\quad {\rm resp.} \\  
&f^{i,{\bf p}}_{\bf a,b} = \Phi^{[i]}(\lambda^{\bf p}_{\bf a,b}) +
[\lambda^{\bf p}_{\bf a,b}]^2 \lambda^{\bf b}_{\bf a,p} 
\lambda^{\bf a}_{\bf b,p} P^{i,{\bf p}}_{\bf a,b}(\lambda^{\bf b}_{\bf a,b},
\lambda^{\bf a}_{\bf b,p})
\end{aligned}$$ 
are independent of the term ${\bf b}$
whenever 
$$\begin{aligned}
&K^{i,{\bf p}}_{\bf c}(t) = t^2(1-t)^2 k^{i,{\bf p}}_{\bf c}(t) ,\\
&\Phi^{[0]}(t) = \Phi(t), \quad 
\Phi^{[1]}(t) = \Phi^{[2]}(t) \equiv \big[ \Psi(t)/(t-1) \big]^\prime
\end{aligned}$$
with arbitrary polynomials $k^{i,{\bf p}}_{\bf c}$ 
and the polynomials $\Phi,\Psi$ have the form $(\theIIIpxxxi.\theiiipxxxi)$\  %$(3.31)$
with arbitrarily fixed polynomials $\Phi_1,\Psi_1$ of one variable.
 
Repeating the calculations of Lemma \theIIIpxix.\theiiipxix,  %$3.19$, 
we see that 
$(\theIIIpxixww.\theiiipxixww)$\   %$(3.19^{\prime\prime})$ 
holds independently of the choice of 
$k^{i,{\bf p}}_{\bf c}, \Phi_1,\Psi_1, R^{i,{\bf p}}_{\bf a,b}$.
Notice that we have constructed the polynomials
$P^{i,{\bf p}}_{\bf a,b}(0,1-t) = P^{i,{\bf p}}_{\bf b,a}(1-t,0)$
in terms of $K^{i,{\bf p}}_{\bf a}$
in a manner such that
$(\theIIIpxxiii.\theiiipxxiii)$\  %$(3.23)$ 
should be fulfilled. 
Thus the expression
$\big[ \overline{\xi}^{\ {\bf b}}_{\bf p,a} \big]^{-1} \!
\left[ x^{[i]}\big( (1 \!-\! t) ({\bf a\!-\! p}) \big) 
[\Phi^{[i]}]^\prime \!(t) 
\xi^{\bf b}_{\bf p,a} +
t^2 (1\!-t) P^{i,{\bf p}}_{\bf a,b}(0,1\!-t) \right]
\big( \!=\! K^{i,{\bf p}}_{\bf a}(t) \big)$
is independent of ${\bf b}$ automatically
which completes the proof in view of Lemma \theIIIpxxiv.\theiiipxxiv.  %3.24. 

\bigskip\bigskip
\noindent\ADS
{\bfs 4. Invariance}

\ADT
\newcounter{IVpi} \setcounter{IVpi}{\value{SEC}}
\newcounter{ivpi} \setcounter{ivpi}{\value{THM}} 
\bigskip
\noindent
{\bf Lemma \NUM.}  %4.1.} 
\it Let ${\bf G}:{\bf x}\mapsto {\bf xA+w}$
be an invertible affine map $\RR^2\leftrightarrow\RR^2$.
A spline procedure 
$\hbox{{\Goth S}}: ({\cal T},F) \mapsto f_{{\cal T},F}$ satisfying Postulate A
is ${\bf G}$-invariant if and only if $(\theIpxi.\theipxi)$ %$(1.11)$ 
holds
for any non-degenerate triangle ${\bf T}$ with distinguished vertex ${\bf p}$.

\bigskip
\noindent
{\bf Proof.} \rm
The ${\bf G}$-invariance of $\hbox{{\Goth S}}$ means that, given any 
triangular mesh ${\cal T}$, the unit functions $f_{{\cal T},F_{i,{\bf p}}}$
$\big(i\!=\!0,\! 1,\! 2;\ {\bf p}\!\in\!{\rm Vert}({\cal T})\big)$ 
corresponding to the gradient data
$F_{i,{\bf p}} := \big\{ ({\bf p},1,{\bf 0})$ if $i\!\!=\!\! 0$, 
$({\bf p},0,{\bf e}^{[i]})$ for $i\!=\!1\! , 2 \big\} \cup 
\big\{ ({\bf q},0,{\bf 0}): 
{\bf p}\!\ne\! {\bf q}\!\in\! {\rm Vert}({\cal T})\big\}$
are transformed by ${\bf G}$ 
as
\ADT
\newcounter{IVpii} \setcounter{IVpii}{\value{SEC}}
\newcounter{ivpii} \setcounter{ivpii}{\value{THM}}
$$\begin{aligned}
&f_{{\cal T},F_{i,{\bf p}}} \circ {\bf G}^{-1} = 
f_{{\bf G}({\cal T}), {\bf G}^\sharp(F_{i,{\bf p}})} \quad (i=0,1,2) \\ 
&{\rm where} \quad 
{\bf G}^\sharp(F_{i,{\bf p}}) \!=\\
&= \big\{ 
\big( {\bf r},[f_{{\cal T},F_{i,{\bf p}}}\!\circ\! {\bf G}^{-1}]({\bf r}),
\nabla [f_{{\cal T},F_{i,{\bf p}}}\!\circ\!{\bf G}^{-1}]({\bf r}) \big): 
{\bf r} \!\in\! {\bf G}({\rm Vert}({\cal T})) \big\} \!= \\ 
&\hskip 7mm  = \big\{ \big( {\bf G(q)}, f_{{\cal T},F_{i,{\bf p}}}({\bf q}),
[\nabla f_{{\cal T},F_{i,{\bf p}}}({\bf q})][{\bf A}^{\rm T}]^{-1} \big) :
{\bf q}\in {\rm Vert}({\cal T}) \big\}  
\end{aligned} 
\eqno(\NUM)$$  %(4.2)$$
with the gradient data 
of the transformed function on the transformed vertices.
Consider any triangle ${\bf T} ={\rm Co}\{ {\bf a,b,p}\} \in{\cal T}$.
Notice that the basic functions over ${\bf T}$ 
are given as restrictions of the 
unit functions. 
In particular 
$f_{{\cal T},F_{0 \!,{\bf p}}} \vert {\bf T} =
\varphi_{\bf p,T}$ and
$f_{{\cal T},F_{j,{\bf p}}} \vert {\bf T} =
\psi^{(j)}_{\bf p,T}$ $(j\!=\!1,2)$.
On the other hand, by Postulate A, for any gradient data
$G$ on ${\rm Vert}\big({\bf G}({\cal T})\big)$ of the transformed mesh,
such that
$\big({\bf G(p)},\omega,[\alpha,\beta]\big),
\big({\bf G(a)},0,{\bf 0}\big),\big({\bf G(b)},0,{\bf 0}\big) \in G,$ 
we have
$f_{{\bf G}({\cal T}),G} = 
\omega\varphi_{\bf G(p),G(T)}+\alpha \psi^{(1)}_{\bf G(p),G(T)}+
\beta \psi^{(2)}_{\bf G(p),G(T)}$.
We can apply this observation to $(\theIVpii.\theivpii)$\  %$(4.2)$ 
with 
$G := {\bf G}^\sharp(F_{i,{\bf p}})$ $(i=0,1,2)$
to conclude that
\ADT
\newcounter{IVpiii} \setcounter{IVpiii}{\value{SEC}}
\newcounter{ivpiii} \setcounter{ivpiii}{\value{THM}}
$$\begin{aligned}
&\varphi_{\bf p,T} \circ {\bf G}^{-1} = 
\varphi_{\bf G(p),G(T)}, \\ 
&\psi^{(j)}_{\bf p,T}\circ {\bf G}^{-1} = 
\alpha_j \psi^{(1)}_{\bf G(p),G(T)} +
\beta_j \psi^{(2)}_{\bf G(p),G(T)} \\ 
& {\rm where} \ \ 
[\alpha_j,\beta_j] = 
[\nabla \phi^{(j)}_{\bf p,T}({\bf p})] [{\bf A}^{\rm T}]^{-1} =
{\bf e}^{[j]} [{\bf A}^{\rm T}]^{-1} \ \  \ (j=1,2) .
\end{aligned}
\eqno(\NUM)$$  %(4.3)$$
Hence the matrix form in $(\theIpxi.\theipxi)$ %(1.11) 
is immediate: $(\theIVpiii.\theivpiii)$\   %$(4.3)$ 
implies that 
$[{\bf A}^{\rm T}]^{-1} = 
\left[ \genfrac{}{}{0pt}{}{\alpha_1 \ \beta_1 }{ \alpha_2 \ \beta_2} \right]$ 
and
$[\psi^{(1)}_{\bf p,T}\circ{\bf G}^{-1},\psi^{(2)}_{\bf p,T}\circ{\bf G}^{-1}] =
[\psi^{(1)}_{\bf G(p),G(T)},\psi^{(2)}_{\bf G(p),G(T)}] {\bf A}^{-1}$.

\ADT
\newcounter{IVpiv} \setcounter{IVpiv}{\value{SEC}}
\newcounter{ivpiv} \setcounter{ivpiv}{\value{THM}}
\bigskip
\noindent
{\bf Corollary \NUM.}   %4.4.} 
\it There is no affine invariant
${\cal C}^1$-spline 
procedure.
satisfying Postulate A.

\bigskip
\noindent{\bf Proof.} \rm
Proceed by contradiction.
Assume the procedure 
$\hbox{{\Goth S}} : ({\cal T},F)\!\mapsto\! f_{{\cal T},F}$
with basic functions $\varphi_{\bf p,T},\psi^{(j)}_{\bf p,T}$
%$\big( {\bf T}$ non-degenerate triangle, ${\bf p}\in{\rm Vert}({\bf T})\big)$
is affine invariant.
Then, in particular, $(\theIpxi.\theipxi)$ %$(1.11)$ 
holds for all transformations
${\bf G}: {\bf x}\mapsto {\bf xA +w}$ with 
${\rm det}({\bf A})\ne 0$ and ${\bf w}\in\RR^2$.
Consider the triangles
%%\ADT
%%\newcounter{IVpv} \setcounter{IVpv}{\value{SEC}}
%%\newcounter{ivpv} \setcounter{ivpv}{\value{THM}}
$${\bf T}_{\bf b} := {\bf G}_{\bf b}({\bf T}) \qquad
\hbox{where \quad ${\bf G}_{\bf b} : {\bf x}\mapsto {\bf xA}_{\bf b}$
with ${\bf A}_{\bf b}:=
\left[ \genfrac{}{}{0pt}{}{1 \ \ \ \ 0 }{ x({\bf b}) \ y({\bf b}) } \right]$.}$$
%%\eqno(\NUM)$$  %(4.5)$$
Then, according to $(\theIpxi.\theipxi)$, %$(1.11)$,
for the points ${\bf b}$ with $y({\bf b})\ne 0$ we have
\ADT
\newcounter{IVpvi} \setcounter{IVpvi}{\value{SEC}}
\newcounter{ivpvi} \setcounter{ivpvi}{\value{THM}} 
$$\big[\psi^{(1)}_{\bf 0,T}\circ{\bf G}_{\bf b}^{-1},
\psi^{(2)}_{\bf 0,T}\circ{\bf G}_{\bf b}^{-1} \big]  {\bf A}_{\bf b} =
\big[ \psi^{(1)}_{{\bf G}_{\bf b}{\bf (0),G}_{\bf b}{\bf (T)}},
\psi^{(2)}_{{\bf G}_{\bf b}{\bf (0),G}_{\bf b}{\bf (T)}} \big] .
\eqno(\NUM )$$  %(4.6)$$ 
Since \ 
${\bf G}_{\bf b}^{-1} : {\bf y}\mapsto {\bf y}{\bf A}_{\bf b}^{-1}$, \ 
in $(\theIVpvi.\theipvi)$\ \   %$(4.6)$ \ 
we can write \  
$\nabla  [\psi^{(j)}_{{\bf 0},{\bf T}} \circ {\bf G}_{\bf b}^{-1}] ({\bf y}) =
\big[\nabla  \psi^{(j)}_{{\bf 0},{\bf T}} ({\bf yA}_{\bf b}^{-1})\big] 
[{\bf A}_{\bf b}^{\rm T}]^{-1}$. 
Therefore, for any ${\bf y}\!\in\!{\bf T}_{\bf b}$ and 
${\bf b}\!\in\!\RR^2$ with $y({\bf b})\!\ne\! 0$,
%%\ADT
%%\newcounter{IVpviw} \setcounter{IVpviw}{\value{SEC}}
%%\newcounter{ivpviw} \setcounter{ivpviw}{\value{THM}}
$$\begin{aligned}
&\nabla \psi^{(1)}_{{\bf 0},{\bf T}_{\bf b}} ({\bf y}) \!=\!
\big[ \nabla  \psi^{(1)}_{{\bf 0},{\bf T}} ({\bf yA}_{\bf b}^{-1}) \big] 
[{\bf A}_{\bf b}^{\rm T}]^{-1} \!+\!
x({\bf b}) 
\big[ \nabla \psi^{(2)}_{{\bf 0},{\bf T}}( {\bf yA}_{\bf b}^{-1}) \big]
[{\bf A}_{\bf b}^{\rm T}]^{-1} ,\\ 
& \nabla \psi^{(2)}_{{\bf 0},{\bf T}_{\bf b}} ({\bf y}) \!=\!
y({\bf b}) \big[ \nabla \psi^{(2)}_{{\bf 0},{\bf T}}( {\bf yA}_{\bf b}^{-1}) \big]
[{\bf A}_{\bf b}^{\rm T}]^{-1} .
\end{aligned}$$
%%\eqno(\NUM)$$  %(4.6^\prime)$$
Observe that the segment ${\rm Co}\{ {\bf 0},{\bf e}^{[1]} \}$ is
a common edge of all the triangles 
${\bf T}_{\bf b}$.  
Hence,
in view of Remark \theIIIpxxii.\theiiipxxii,  %3.11,  
the gradients
$\nabla \psi^{(j)}_{{\bf 0},{\bf T}_{\bf b}}({\bf y}_t)$ 
with ${\bf y}_t:=t {\bf e}^{[1]}$
must be independent of ${\bf b}$ for any fixed $t\in[0,1]$.
Since \ ${\bf y}_t {\bf A}_{\bf b}^{-1} = {\bf y}_t$ \
$(t\!\in\!\RR,\ y({\bf b})\!\ne\! 0)$, 
our indirect assumption leads to the conclusions that \ 
${\bf 0} = \nabla \psi^{(2)}_{{\bf 0},{\bf T}_{\bf b}} ({\bf y}_t)=
\nabla \psi^{(2)}_{{\bf 0},{\bf T}}( {\bf y}_t )$
and 
$\big[ \nabla \psi^{(2)}_{{\bf 0},{\bf T}_{\bf b}} ({\bf y}_t) \big] 
{\bf A}_{\bf b}^{\rm T} =
\nabla \psi^{(2)}_{{\bf 0},{\bf T}}( {\bf y}_t )$ 
for all $t\!\in\![0,1]$ and ${\bf b}\in \RR^2$ with $y({\bf b})\!\ne\! 0$.
This latter identity means in particular that
$x({\bf b}) \frac{\partial }{ \partial x} 
\psi^{(2)}_{{\bf 0},{\bf T}_{\bf b}}( {\bf y}_t) +
y({\bf b}) \frac{\partial}{\partial y} 
\psi^{(2)}_{{\bf 0},{\bf T}_{\bf b}}( {\bf y}_t) =
\frac{\partial}{\partial y} \psi^{(2)}_{{\bf 0},{\bf T}}( {\bf y}_t)$ 
which is possible with ${\bf b}$-independent 
$\nabla \psi^{(2)}_{{\bf 0},{\bf T}_{\bf b}}( {\bf y}_t)$
only if
$\frac{\partial}{\partial x} 
\psi^{(2)}_{{\bf 0},{\bf T}_{\bf b}}( {\bf y}_t) =
\frac{\partial}{\partial y} 
\psi^{(2)}_{{\bf 0},{\bf T}_{\bf b}}( {\bf y}_t) =0$
$(t\!\in\! [0,1])$.
However, hence we get
${\bf 0}= \nabla \psi^{(2)}_{{\bf 0},{\bf T}_{\bf b}}( {\bf y}_0)$
which  contradicts the defining relations $(\theIpv.\theipv)$ %$(1.5)$ 
with
$\nabla \psi^{(2)}_{{\bf 0},{\bf T}_{\bf b}}( {\bf 0}) =[0,1]$.

\ADT
\newcounter{IVpvii} \setcounter{IVpvii}{\value{SEC}}
\newcounter{ivpvii} \setcounter{ivpvii}{\value{THM}}
\bigskip
\noindent
{\bf Lemma \NUM.}  %4.7.} 
{\rm (Reflection lemma).}
\it  Let ${\bf T}$ be a non-degenarate triangle of the form
${\bf T}:= {\rm Co} \{ {\bf 0},\rho {\bf e}^{[1]},{\bf b} \}$
and assume $\hbox{\Goth S}$ is a spline procedure satisfying Postulate A. 
Then, for the fixed points ${\bf u}_t := t{\bf e}^{[1]}$
of the reflection 
${\bf K}={\bf K}^{-1} : {\bf x}\mapsto [x({\bf x},-y({\bf x}] ={\bf xU}$,
${\bf U}=\big[ \genfrac{}{}{0pt}{}{ 1 \ \ \ 0 }{ 0 \ -1}\big]$
through the $x$-axis we have
%%\ADT
%%\newcounter{IVpviii} \setcounter{IVpviii}{\value{SEC}}
%%\newcounter{ivpviii} \setcounter{ivpviii}{\value{THM}} 
$$\big\langle \nabla \varphi_{\bf 0,T}({\bf u}_t) 
\big\vert {\bf e}^{[2]} \big\rangle \!=\! 0, \
\big\langle \nabla \psi^{(1)}_{\bf 0,T}({\bf u}_t) 
\big\vert {\bf e}^{[2]} \big\rangle \!=\! 0, \  
\psi^{(2)}_{\bf 0,T}({\bf u}_t) \!=\! 0 
\ \ (t\!\in\![0,\rho]).$$
%%\eqno(\NUM)$$  %(4.8)$$
{\bf Proof.} \rm The triangles ${\bf T}$ and ${\bf K(T)}$ 
are adjacent,
the segment ${\rm Co}\{ {\bf 0},\rho{\bf e}^{[1]}\}$
is their common edge. 
According to Remark \theIIIpxxii.\theiiipxxii,  %3.11, 
the pairs
$\varphi_{\bf 0,T},\varphi_{\bf 0,K(T)}$ resp.
$\psi^{(j)}_{\bf 0,T},\psi^{(j)}_{\bf 0,K(T)}$ 
of basic functions
must be coupled ${\cal C}^1$-smoothly along it:
$\varphi_{\bf 0,T}({\bf u}_t) \!\!=\! \varphi_{\bf 0,K(T)}({\bf u}_t)$,
$\nabla\varphi_{\bf 0,T}({\bf u}_t) \!=\! 
\nabla\varphi_{\bf 0,K(T)}({\bf u}_t)$
resp.
$\psi^{(j)}_{\bf 0,T}({\bf u}_t) \!=\! \psi^{(j)}_{\bf 0,K(T)}({\bf u}_t)$,
$\nabla\psi^{(j)}_{\bf 0,T}({\bf u}_t) = 
\nabla\psi^{(j)}_{\bf 0,K(T)}({\bf u}_t)$
for all $t\in[0,\rho]$.
On the other hand, the transformation rules $(\theIpxi.\theipxi)$ %$(1.11)$ 
require
$\varphi_{\bf 0,K(T)} \!=\! \varphi_{\bf 0,T}\circ {\bf K}^{-1}$ resp.
$\big[ \psi^{(1)}_{\bf 0,K(T)}, \psi^{(2)}_{\bf 0,K(T)} \big] \!=\!
\big[ \psi^{(1)}_{\bf 0,T}\circ {\bf K}^{-1}, 
\psi^{(2)}_{\bf 0,T}\circ {\bf K}^{-1} \big]{\bf U}$ \ i.e. \ 
$\varphi_{\bf 0,K(T)}({\bf y}) = \varphi_{\bf 0,T}({\bf yU})$ and
$\psi^{(j)}_{\bf 0,K(T)}({\bf y})= (-1)^{j-1}\psi^{(j)}_{\bf 0,T}({\bf y})$
for all ${\bf y}\in{\bf K(T)=TU}$.  
By passing to gradients, since ${\bf K}={\bf K}^{-1}$ and 
${\bf U}={\bf U}^{-1} ={\bf U}^{\rm T}$, we get \
$\nabla \varphi_{\bf 0,K(T)} ({\bf y}) = 
\nabla [\varphi_{\bf 0,T}({\bf yU})]{\bf U}$ \ 
and \ 
$\nabla \psi^{(j)}_{\bf 0,K(T)}({\bf y}) = 
(-1)^{j-1}[\nabla \psi^{(j)}_{\bf 0,T}({\bf y})$ \
for the points \  ${\bf y}\in{\bf K(T)}$.
In particular on the common edge of \ ${\bf T}$ with ${\bf K(T)}$ \  
we must have
$\nabla \varphi_{\bf 0,T} ({\bf u}_t)
\nabla \varphi_{\bf 0,K(T)} ({\bf u}_t) = 
\nabla [\varphi_{\bf 0,T}(({\bf u}_t))]{\bf U}$
i.e. 
$\frac{\partial}{ \partial x^{[k]}} \varphi_{\bf 0,T}({\bf u}_t) 
= (-1)^{k-1} \frac{\partial}{ \partial x^{[k]}} \varphi_{\bf 0,T}({\bf u}_t)$
implying 
$0 = \frac{\partial}{ \partial y} \varphi_{\bf 0,T}({\bf u}_t) =
\big\langle \nabla \varphi_{\bf 0,T}({\bf u}_t) 
\big\vert {\bf e}^{[2]} \big\rangle$.
Similarly we conclude that
$\psi^{(j)}_{\bf 0,T}({\bf u}_t) = 
\psi^{(j)}_{\bf 0,K(T)}({\bf u}_t) =
(-1)^{j-1} \psi^{(j)}_{\bf 0,T}({\bf u}_t)$
implying $\psi^{(2)}_{\bf 0,K(T)}({\bf u}_t) =0$
and
$\nabla \psi^{(j)}_{\bf 0,T}({\bf u}_t) \!=\!
\nabla \psi^{(j)}_{\bf 0,K(T)}({\bf u}_t) \!=\! 
(-1)^{j-1} [\nabla \psi^{(j)}_{\bf 0,T}({\bf u}_t)]{\bf U}$
implying in particular
$\frac{\partial}{ \partial y} \psi^{(1)}_{\bf 0,T}({\bf u}_t) \!=\! 0$. 

\ADT
\newcounter{IVpix} \setcounter{IVpix}{\value{SEC}}
\newcounter{ivpix} \setcounter{ivpix}{\value{THM}}
\bigskip
\noindent
{\bf Proposition \NUM.}  %4.9.} 
\it Homothetically invariant ${\cal C}^1$-spline
procedures satisfying Postulate A are shape uniform on edges
$($i.e. they satisfy Postulate B automatically$)$.

\bigskip
\noindent
{\bf Proof.} \rm
Let ${\bf T}:={\rm Co}\{ {\bf 0},{\bf e}^{[1]},{\bf e}^{[2]}\}$ 
and define
\ADT
\newcounter{IVpxviii} \setcounter{IVpxviii}{\value{SEC}}
\newcounter{ivpxviii} \setcounter{ivpxviii}{\value{THM}} 
$$\Phi(t):= \varphi_{\bf 0,T}({\bf u}_t), \quad 
\Psi(t) := \phi^{(1)}_{\bf 0,T}({\bf u}_t) \quad 
{\rm where} \ \
{\bf u}_t:= (1\!-\! t) {\bf e}^{[1]} .% \ t\!\in\! [0,1].
\eqno(\NUM)$$    %(4.18)$$
Cosider any other non-degenerate triangle
$\widetilde{\bf T} := {\rm Co}\{ {\bf p,a,b}\}$.
Due to the arbitrariness of the choice of $\widetilde{\bf T}$,
it suffices to see only that,
for $j \!=\!1,2$ and $t\!\in\![0,1]$,
$$\varphi_{{\bf 0},\widetilde{\bf T}} ({\bf y}_t) \!=\! \Phi(t), \ \ 
\psi^{(j)}_{{\bf 0},\widetilde{\bf T}} ({\bf y}_t) \!=\!
{\rm Const}^{(j)}_{\widetilde{\bf T}}\Psi(t)
\quad {\rm with} \quad  
{\bf y}_t \!:=\! (1\!-\! t){\bf a} \!+\! t{\bf p} .
\eqno(4.9)$$
It is a crucial fact that we can find a homothetic transformation
$${\bf G} :{\bf x}\mapsto {\bf xA}+{\bf p} \quad\hbox{such that}\ 
{\bf G}({\bf e}^{[1]}) ={\bf a}, \ \ 
{\bf G}({\bf T}) \cap \widetilde{\bf T} ={\rm Co}\{ {\bf a,p}\} .
\eqno(4.10)$$
Actually ${\bf A} = \big[ 
\genfrac{}{}{0pt}{}{ x({\bf a-p}) \ \ \ y({\bf a-p}) }{ 
-\sigma y({\bf a-p}) \ \sigma x({\bf a-p}) }
\big]$ where $\sigma=1$ if the points ${\bf a,b,p}$ are oriented
clockwise and $\sigma=-1$ else. 
According to $(\theIpxi.\theipxi)$, %$(1.11)$, 
$\varphi_{\bf G(0),G(T)} = \varphi_{\bf 0,T}\circ {\bf G}^{-1}$.
Since ${\rm Co}\{ {\bf p,a}\}$ is a common edge of 
${\bf T}$ and $\widetilde{\bf T}$, 
in view of Remark \theIIIpxxii.\theiiipxxii\  %3.11 
we have
\ADT
\newcounter{IVpxi} \setcounter{IVpxi}{\value{SEC}}
\newcounter{ivpxi} \setcounter{ivpxi}{\value{THM}}
$$\varphi_{{\bf p},\widetilde{\bf T}} ({\bf y}_t) =
\varphi_{\bf p,G(T)}({\bf y}_t) = 
\varphi_{\bf 0,T}\big({\bf G}^{-1}({\bf y}_t)\big) =
\varphi_{\bf 0,T}\big({\bf u}_t\big) = \Phi(t) 
\eqno(\NUM)$$  %(4.11)$$
which proves the first part of $(\theIIIpx.\theiiipx)$.  %$(4.10)$. 
To prove 
$\psi^{(j)}_{bf 0,\widetilde{T}}({\bf y}_t) =
{\rm Const}^{(j)}_{\widetilde{\bf T}}\Psi(t)$,
consider also the symmetry
\ADT
\newcounter{IVpxii} \setcounter{IVpxii}{\value{SEC}}
\newcounter{ivpxii} \setcounter{ivpxii}{\value{THM}}
$$\begin{aligned}
&{\bf H} %\!=\! {\bf H}^{-1} \!\! 
:{\bf x} \!\mapsto\! 
[y({\bf x}),x({\bf x})]={\bf xS}, 
\ \ {\bf v}_t:= {\bf u}_t{\bf S} =t{\bf e}^{[2]} 
\ \hbox{where \ 
${\bf S} \!:=\! \big[ \genfrac{}{}{0pt}{}{0 \ 1 }{ 1 \ 0 }\big]$} 
\end{aligned}
\eqno(\NUM)$$  %(4.12)$$
of the triangle ${\bf T}$. 
By $(\theIpxi.\theipxi)$ %$(1.11)$  
we have \  
$\big[ \psi^{(1)}_{\bf 0,T}, \psi^{(2)}_{\bf 0,T} \big] \!=\!\!
\big[ \psi^{(1)}_{\bf 0,T} \circ {\bf H}, 
\psi^{(2)}_{\bf 0,T}  \circ {\bf H} \big] {\bf S}$ 
whence 
$$\big[ \psi^{(1)}_{\bf 0,T} ({\bf u}_t), \psi^{(2)}_{\bf 0,T} ({\bf u}_t)\big] \! =
\big[ \psi^{(1)}_{\bf 0,T} ({\bf v}_t) , 
\psi^{(2)}_{\bf 0,T}  ({\bf v}_t) \big] 
\big[ \genfrac{}{}{0pt}{}{ 0 \ 1 }{ 1 \ 0} \big] = 
\big[ \psi^{(2)}_{\bf 0,T} ({\bf v}_t) , 
\psi^{(1)}_{\bf 0,T}  ({\bf v}_t) \big] .$$
Thus 
$\psi^{(2)}_{\bf 0,T} ({\bf v}_t) = \Psi(t)$
while 
$\psi^{1)}_{\bf 0,T} ({\bf v}_t) = \psi^{(2)}_{\bf 0,T} ({\bf u}_t)$. 
On the other hand, by Lemma \theIVpvii.ivpvii,   %4.7, 
$\psi^{(2)}_{\bf 0,T} ({\bf u}_t) =0$. 
Finally we apply Remark \theIIIpxxii.\theiiipxxii\   %3.11 
and $(\theIpxi.\theipxi)$ %$(1.11)$ 
to the points ${\bf y}_t$ of the common edge 
${\rm Co} \{ {\bf a,p}\}$ between the triangles
${\bf G(T)}$ and $\widetilde{\bf T}$. 
Hence we conclude that
\ADT
\newcounter{IVpxiii} \setcounter{IVpxiii}{\value{SEC}}
\newcounter{ivpxiii} \setcounter{ivpxiii}{\value{THM}}
$$\begin{aligned}
&\big[ \psi^{(1)}_{{\bf p},\widetilde{\bf T}}({\bf y}_t) ,
\psi^{(2)}_{{\bf p},\widetilde{\bf T}}({\bf y}_t) \big] =
\big[ \psi^{(1)}_{{\bf G(0)},{\bf G(T)}}({\bf y}_t) ,
\psi^{(2)}_{{\bf G(0)},{\bf G(T)}}({\bf y}_t) \big] = \cr 
&=\! \big[ \psi^{(1)}_{{\bf 0},{\bf T}}({\bf u}_t) ,
\big[ \psi^{(2)}_{{\bf 0},{\bf T}}({\bf u}_t) \big] {\bf A} \!=\!
\big[ \Psi(t),0\big] {\bf A} \!=\! 
\Psi(t) \big[ x({\bf a \!-\! p} , y({\bf a \!-\! p} \big] .
\end{aligned}
\eqno(\NUM)$$   %(4.13)$$
Thus \  
$\psi^{(j)}_{{\bf p},\widetilde{\bf T}}({\bf y}_t) =
x^{[j]}({\bf a-p}) \Psi(t)$ $(j\!=\! 1,2)$ \  
which completes the proof.

\ADT
\newcounter{IVpxiv} \setcounter{IVpxiv}{\value{SEC}}
\newcounter{ivpxiv} \setcounter{ivpxiv}{\value{THM}}
\bigskip
\noindent
{\bf \NUM.  %4.14. 
Proof of Therem \theIIpxi.\theiipxi}    %2.11}

\medskip
It is clear that the coordinate values 
$\xi^{\bf v}_{\bf p,a} = 
\langle {\bf v \!-\! a} \vert {\bf p\!-\! a} \rangle /
\Vert {\bf p\!-\! a}\Vert^2$ are homothetic invariant
i.e. 
$\xi^{\bf G(v)}_{\bf G(p),G(a)} = \xi^{\bf v}_{\bf p,a}$ 
whenever the transformation
${\bf G}:\RR^2\leftrightarrow\RR^2$
is of the form 
${\bf G(x) = w + \rho (x-q)S}$ with a constant $\rho>0$ and an
orthogonal $2\times 2$-matrix ${\bf S}$.
The baricentric coordinates $\lambda^{\bf p}_{\bf T}$ are even
affine invariant as it is well-known from classical Projective Geometry.
Hence it suffices to see that the invariance relations 
$$\varphi_{\bf G(p),G(T)} = \varphi\circ {\bf G}^{-1}, \ \
\big[ \psi^{(1)}_{\bf G(p),G(T)}, \psi^{(2)}_{\bf G(p),G(T)} \big] =
\big[ \psi^{(1)}\circ {\bf G}^{-1}, \psi^{(1)}\circ {\bf G}^{-1} \big]
{\bf S}$$
imply that 
$k^{i,{\bf p}}_{\bf a} \equiv 0$ $(i=0,1,2)$
whenever 
${\bf T} = {\rm Co}\{ {\bf a,b,p} \}$ 
is a non-degenerate triangle
and ${\bf G}:\RR^2\leftrightarrow\RR^2$ is 
the orthogonal reflection through ${\rm Line}\{ {\bf a,p}\}$
i.e.
$${\bf G(a)=a, \ \ G(p)=p}, \quad 
\xi^{\bf G(x)}_{\bf p,a} = \xi^{\bf x}_{\bf p,a}, \ \ 
\overline{\xi}^{\bf G(x)}_{\bf p,a} = 
- \overline{\xi}^{\bf x}_{\bf p,a} \quad ({\bf x}\in\RR^2) ,$$
so that   
${\bf G(x) = a \!+\! (x \!-\! a)S}$  where 
${\bf S} = \Vert {\bf p\!-\! a}\Vert^{-2}$
\DEL{\left[ 
\begin{matrix} 
{\bf p\!-\! a}\\ {\bf (p\!-\! a)R} \end{matrix} \right]^{\rm T}
\left[ \begin{matrix} 1 &0 \\ 0 &-1 \end{matrix} \right]
\left[ \begin{matrix} {\bf p\!-\! a}\\ {\bf (p\!-\! a)R} \end{matrix} \right].
}%ENDDEL
$\genfrac[]{0pt}{1}{{\bf p\!-\! a}}{ {\bf (p\!-\! a)R}}^{\rm T}
\genfrac[]{0pt}{1}{{0\ 1}}{{1\ 0}}
\genfrac[]{0pt}{1}{{\bf p\!-\! a}}{ {\bf (p\!-\! a)R}}$.
Let us first investigate the relation
$\varphi_{\bf G(p),G(T)} = \varphi_{\bf p,T}\circ {\bf G}^{-1}$.
By pluging 
$$\begin{aligned}
&t=t({\bf x}) := \lambda^{\bf a}_{\bf G(T)}({\bf x}) =
\lambda^{\bf G(a)}_{\bf G(T)}({\bf x}) =
\lambda^{\bf a}_{\bf T}\big({\bf {\bf G}^{-1}(x)}\big),\\ 
&s=s({\bf x}) := \lambda^{\bf G(b)}_{\bf G(T)}({\bf x}) =
\lambda^{\bf b}_{\bf T}\big({\bf {\bf G}^{-1}(x)}\big)
\end{aligned}$$
in the expressions of 
$\varphi_{\bf G(p),G(T)}$ resp. $\varphi_{\bf p,T}\circ {\bf G}^{-1}$
formed with $(\theIIpviii.\theiipviii)$,  %$(2.8)$, 
since 
$\lambda^{\bf p}_{\bf T} = 
1 - \lambda^{\bf a}_{\bf T} -\lambda^{\bf b}_{\bf T}$
and ${\bf G(a)=a}$ resp. ${\bf G(p)=p}$, 
we get
$$\begin{aligned}
&0 = \varphi_{\bf G(p),G(T)}({\bf x})  - 
\varphi_{\bf p,T}\circ {\bf G}^{-1} ({\bf x}) =
\varphi_{\bf p,G(T)}({\bf x})  - 
\varphi_{\bf p,T}\circ {\bf G}^{-1} ({\bf x}) =\\ 
&\ \, = \Big[ \Phi(1 \!-\! s \!-\! t) \!+\! 
(1 \!-\! s \!-\! t)^2 st P^{\bf p}_{\bf a,G(b)} (s,t) \Big] \!-\\ 
&\hskip20mm - \Big[ \Phi(1 \!-\! s \!-\! t) \!+\! 
(1 \!-\! s \!-\! t)^2 st P^{\bf p}_{\bf a,b}(s,t)\Big],\\ 
&0 =  
P^{\bf p}_{\bf a,G(b)} (s,t) - 
P^{\bf p}_{\bf a,b} (s,t) = 
\left[ 
s \Big\{ \xi^{\bf a}_{\bf p, G\!(b)} \frac{\Phi^\prime(1 \!-\! s)
}{ s(1\!-\! s)^2}\!+\!
\overline{\xi}^{\bf a}_{\bf p,G\!(b)} k^{0,{\bf p}}_{\bf G(b)}(s) \!\Big\} 
\right. +\\
&\hskip 20mm + \left. 
t \Big\{ \xi^{\bf G\!(b)}_{\bf p,a}
\frac{\Phi^\prime(1 \!-\! t)}{ t(1 \!-\! t)^2} \!+\!
\overline{\xi}^{\bf G\!(b)}_{\bf p,a} k^{0,{\bf p}}_{\bf a}(t) \!\Big\}
\!+\! st R^{0,{\bf p}}_{\bf a,G\!(b)} (s,t) \right] \!\!-\! \\ 
&\hskip 20mm - 
\left[ 
s \Big\{ \xi^{\bf a}_{\bf p,b}
\frac{\Phi^\prime(1 \!-\! s) }{ s(1\!-\! s)^2}\!+\!
\overline{\xi}^{\bf a}_{\bf p,b} k^{0,{\bf p}}_{\bf b}(s) \!\Big\} 
\right.\!+\\
&\hskip20mm + \left. t \Big\{ \xi^{\bf b}_{\bf p,a}
\frac{\Phi^\prime(1 \!-\! t) }{ t(1 \!-\! t)^2} \!+\!
\overline{\xi}^{\bf b}_{\bf p,a} k^{0,{\bf p}}_{\bf a}(t) \Big\}
\!+\! st R^{0,{\bf p}}_{\bf a,b} (s,t) \right] 
\end{aligned}$$
Comparing the coefficients of the monomials
$s^m t^n$, in view of $(2.9)$ we see 
that
$$\begin{aligned}
{\rm (i)} \ 
0 &\!=\! \Big[ \xi^{\bf G\!(b)}_{\bf p,a} \frac{\Phi^\prime(1-t) }{ t(1-t)^2} +
\overline{\xi}^{\bf G\!(b)}_{\bf p,a} k^{0,{\bf p}}_{\bf a}(t) \Big] -
\Big[ \xi^{\bf b}_{\bf p,a} \frac{\Phi^\prime(1-t) }{ t(1-t)^2} +
\overline{\xi}^{\bf b}_{\bf p,a} k^{0,{\bf p}}_{\bf a}(t) \Big] ,\\ 
{\rm (ii)} \ 
0 &\!=\! \Big[ \xi^{\bf a}_{\bf p,G\!(b)} 
\frac{\Phi^\prime(1-s) }{ s(1 \!-\! s)^2} \!+\!
\overline{\xi}^{\bf a}_{\bf p,G\!(b)} k^{0,{\bf p}}_{\bf G\!(b)}(s) \Big] \!-\!
\Big[ \xi^{\bf a}_{\bf p,b} \frac{\Phi^\prime(1-s) }{ s(1 \!-\! s)^2} \!+\!
\overline{\xi}^{\bf a}_{\bf p,b} k^{0,{\bf p}}_{\bf a}(s) \Big] ,\\
{\rm (iii)} \
0 &\!=\! R^{0,{\bf p}}_{\bf a,G\!(b)} (s,t) - R^{0,{\bf p}}_{\bf a,b} (s,t) .
\end{aligned}$$
Since 
$\xi^{\bf G\!(b)}_{\bf p,a} = \xi^{\bf b}_{\bf p,a}$ and
$\overline{\xi}^{\bf G\!(b)}_{\bf p,a} = -\overline{\xi}^{\bf b}_{\bf p,a}$,
from (i) we conclude that $k^{0,{\bf p}}_{\bf a}(t) =0$.
On the other hand, (iii) implies the isometry invariance of the
non-principal parts because 
the lines ${\rm Line}\{{\bf p,a)}\}$ 
can be chosen arbitrarily and hence the corresponding reflections
generate the whole group of self-isometries of $\RR^2$.

The treatment of the relations
$$\big[ \psi^{(1)}_{\bf G(p),G(T)}, \psi^{(2)}_{\bf G(p),G(T)} \big] =
\big[ \psi^{(1)}\circ {\bf G}^{-1}, \psi^{(1)}\circ {\bf G}^{-1} \big]
{\bf S}$$ 
is analogous by using the vectorial forms
$\ppsi_{\bf q,W} := \big[ \psi^{(1)}_{\bf q,W}, \psi^{(2)}_{\bf q,W} \big]$
for triangles ${\bf W}$ with distinguished vertex ${\bf q}$.
With this formalism the above invariance relation can be written as \ 
$\ppsi_{\bf G\!(p),G\!(T)} ({\bf x}) = 
\big[ \ppsi_{\bf p,T}\big( {\bf G}^{-1}({\bf x})\big) \big]{\bf S}$ \ 
where ${\bf G}={\bf G}^{-1}$, ${\bf G\!(p)=p,\ G\!(a)=a}$ and
$$\begin{aligned}
&\ppsi_{\bf G\!(p),G\!(T)} ({\bf x}) = 
\Theta(1\!-\! s\!-\! t) ({\bf x \!-\! p}) +
(1\!-\! s\!-\! t)^2 s t Q^{\bf p}_{\bf a,G\!(b)}(s,t) ,\\
&\big[ \ppsi_{\bf p,T}\big( {\bf G}^{-1}({\bf x})\big) \big]{\bf S} =
\Theta(1\!-\! s\!-\! t)\big( {\bf G\!(x)} \!-\! {\bf p}\big){\bf S} +
(1\!-\! s\!-\! t)^2 s t Q^{\bf p}_{\bf a,b}(s,t) {\bf S}  
\end{aligned}$$
with the vector valued polynomials
$$\begin{aligned}
\QQ^{\bf p}_{\bf a,w} (s,t) &:= 
\big[ Q^{1,{\bf p}}_{\bf a,w} (s,t) , Q^{2,{\bf p}}_{\bf a,w} (s,t) \big] =\\
&= s \Big\{ \xi^{\bf a}_{\bf p,w}
\frac{\Theta^\prime(1 \!-\! s) }{ s(1\!-\! s)^2} {\bf (w \!-\! p)} + 
\overline{\xi}^{\bf a}_{\bf p,w}  
{\bf k}_{\bf w}^{{\bf p}}(s) \Big\} +\\  
&+ \ \ t \Big\{ \xi^{\bf w}_{\bf p,a}  
\frac{\Theta^\prime(1 \!-\!t) }{ t(1 \!-\! t)^2} {\bf (a \!-\! p)} +
\overline{\xi}^{\bf w}_{\bf p,a}
{\bf k}_{\bf a}^{{\bf p}}(t) \Big\} 
+  st {\bf R}^{{\bf p}}_{\bf a,w}(s,t) 
\end{aligned}$$
for \ ${\bf w}:= {\bf b,G\!(b)}$ \ 
where \ ${\bf R}^{{\bf p}}_{\bf a,w} := 
\big[ {\bf R}^{1,{\bf p}}_{\bf a,w}, {\bf R}^{2,{\bf p}}_{\bf a,w} \big]$ \ 
and \ ${\bf k}_{\bf u}^{{\bf p}} :=
\big[ {\bf k}_{\bf w}^{1,{\bf p}},{\bf k}_{\bf w}^{2,{\bf p}} \big]$.
Clearly  
${\bf (G\!(x) \!- \! p)S} \!=\! {\bf (G\!(x) \!-\! G\!(p))S} \!=\!
{\bf \big( (a + (x\!-\! a)S) \!-\! (a + (p\!-\! a)S) \big)S} \!=\!
{\bf (x \!-\! p)S^2} = {\bf x\!-\! p}$. 
Hence 
the comparison of the coefficients of the monomials $s^m t^n$
yields
$$\begin{aligned}
&{\rm (i^\prime)} \quad  
0 = \Big[ \xi^{\bf G\!(b)}_{\bf p,a} \frac{\Theta^\prime(1-t) }{ t(1-t)^2}
{\bf \big( G\!(b) \!-\! p\big)} \!+\!
\overline{\xi}^{\bf G\!(b)}_{\bf p,a} {\bf k}^{{\bf p}}_{\bf a}(t) \Big] \!-\\
&\hskip 20mm -
\Big[ \xi^{\bf b}_{\bf p,a} \frac{\Theta^\prime(1-t) }{ t(1-t)^2} 
{\bf (b \!-\!p)} \!+\!
\overline{\xi}^{\bf b}_{\bf p,a} {\bf k}^{{\bf p}}_{\bf a}(t) \Big]{\bf S} ,
\\ 
&{\rm (ii^\prime)} \quad 
0 = \Big[ \xi^{\bf a}_{\bf p,G\!(b)} \frac{\Theta^\prime(1 \!-\!s) 
}{ s(1 \!-\! s)^2} {\bf (a\!-\! p)} \!+\!
\overline{\xi}^{\bf a}_{\bf p,G\!(b)} {\bf k}^{{\bf p}}_{\bf G\!(b)}(s) \Big] 
\!-\\ 
&\hskip 20mm - \Big[ \xi^{\bf a}_{\bf p,b} \frac{\Theta^\prime(1 \!-\!s)}{ s(1\!-\! s)^2} 
{\bf (a\!-\! p)} \!+\!
\overline{\xi}^{\bf a}_{\bf p,b} {\bf k}^{{\bf p}}_{\bf a}(s) \Big]{\bf S} ,
\\
&{\rm (iii^\prime)} \quad
0 = {\bf R}\strut^{{\bf p}}_{\bf a,G\!(b)} (s,t) - 
{\bf R}^{{\bf p}}_{\bf a,b} (s,t) {\bf S} .
\end{aligned}$$
Considering again ${\rm (i^\prime)}$, since
${\bf G\!(b) \!-\! p} = {\bf G\!(b) \!-\! G\!(p)} =
{\bf (b \!-\! p)}S$ 
and since 
$\xi^{\bf G\!(b)}_{\bf p,a} = \xi^{\bf b}_{\bf p,a}$ resp.
$\overline{\xi}^{\bf G\!(b)}_{\bf p,a} = -\overline{\xi}^{\bf b}_{\bf p,a}$, 
we conclude ${\bf k}^{\bf p}_{\bf a}(t)=0$.

\bigskip
\noindent
{\bf Acknowledgement.} 
This research was supported by the Ministry of Human Capacities,
Hungary grant 20391-3/2018/FEKUSTRAT.
     
\bigskip
\bigskip
\noindent
{\bfs REFERENCES}

\bigskip

\begin{itemize}
\item[[\Barn]]
R.E. Barnhill, Smooth interpolation over triangles (p.47) 45-70?, in:
R.E. Barnhild and R.F. Riesenfeld (Ed.-s),
Computer Aided Geometric Design, Academic Press, NY-San Francisco-London,1974.

%\smallskip
\item[[\Cox]] 
H-S.M. Coxeter, Introduction to Geometry, New York: Wiley, 1969. 

%\smallskip
\item[[\Fulton]]
William Fulton, Algebraic Curves, 
Mathematics Lecture Note Series. W.A. Benjamin. p. 112. (1974).
%ISBN 0-8053-3081-4.

%\smallskip
\item[[\Hamilton]]
H. Hamilton, Plane Algebraic Curves, Calderon Press, Oxford, 1920. 

%\smallskip
\item[[\Lorenz]] 
Lorentz, Hermite interpolation by algebraic
polynomials: A survey in Numerical Analysis 2000 (Elsevier).

%\smallskip
\item[[\Meng]]
Meng wu, Bernard Mourrain, André Galligo, B. Nkonga, 
Hermite Type Spline Spaces over Rectangular Meshes with Complex 
Topological Structures,
%March 2017
Communications in Computational Physics 21/3, 835-866.
%DOI: 10.4208/cicp.OA-2016-0030

%\smallskip
\item[[\Schumacher]] 
%https://epubs.siam.org/doi/book/10.1137/1.9781611973907 
L.L. Schumacher, Spline Functions, Computational Methods, SIAM 2015. 

%\smallskip
\item[[\Serg]]
I.V. Sergienko, O.M. Lytvyn, O.O. Lytvyn and O.I. Denisova,
Explicit formulas for interpolating splines of degree 5 on the triangle,
Cybernetics and Systems Analysis, Vol. 50, No.5 (2014), 670-678.
%DOI 10.1007/s10559-014-9657-x

%\smallskip
\item[[\Zlam]] 
M. Zl\'amal and A. \v Zeni\v sek, Mathematical aspects of the FEM, in
Technical Physical and Mathematical Priciples of the FEM,
Ed.-s V. Kola\v z et al., Akademia, Praha, 1971, pp. 15-39.

%\smallskip
%\item[[\St]] Stach\'o et al., Triang.mws,
%Home page http://www.math.u-szeged.hu/~stacho/.

\end{itemize}

\bigskip
\bigskip
L.L. STACH\'O

Bolyai Institute, 

Interdisciplinary Excellence Centre, 

University of Szeged

{\tt stacho@math.u-szeged.hu}

\end{document}